\documentclass[11pt]{article}
\usepackage{latexsym}
\usepackage{amsmath,amsthm,lipsum}
\usepackage{amssymb,amsfonts,amscd,mathtools}
\usepackage{graphicx}
\usepackage{tikz,tikz-cd}
\usetikzlibrary{tikzmark}
\usepackage{multirow}
\usepackage{mathabx} 
\usepackage[dvipsnames]{xcolor}
\usepackage{hyperref}
\hypersetup{
	colorlinks,
	linkcolor={MidnightBlue},
	citecolor={Green}
}
\usepackage[nameinlink]{cleveref}
\Crefname{algocf}{Algorithm}{Algorithms}
\usepackage[footnotesize]{caption}
\usepackage{subcaption}
\usepackage{cite}
\usepackage[ruled,noend]{algorithm2e}
\numberwithin{equation}{section}
\usepackage[textsize=footnotesize, backgroundcolor=gray!25]{todonotes}
\usepackage{changes}
\usepackage{hhline}
\usepackage{tabularx}
\usepackage{array}
\usepackage{siunitx}
\usepackage{enumitem}

\usepackage{soul}
\usepackage{thmtools} 

\usepackage{diagbox}

\theoremstyle{plain}
    \newtheorem{theorem}{Theorem}[section]
    \newtheorem{proposition}[theorem]{Proposition}
    \newtheorem{corollary}[theorem]{Corollary}

\theoremstyle{definition}
    \newtheorem{definition}[theorem]{Definition}

\newtheorem{example}[theorem]{Example}

\AtBeginEnvironment{example}{%
  
  \pushQED{\qed}
}

\AtEndEnvironment{example}{%
  \popQED
}

\theoremstyle{remark} 
    \newtheorem{remark}[theorem]{Remark}

\newcommand{\sH}{{\mathcal H}}

\newcommand{\sW}{{\mathcal W}}

\newcommand\bC{{\mathbb C}}

\newcommand\bR{{\mathbb R}}

\newcommand\suchthat{ \mid }

\DeclareMathOperator*{\Hess}{Hess}

\begin{document}
\title{Numerical Elimination: Computing Complements\\ of Real Hypersurfaces Using Pseudo-Witness Sets}

\author{
Paul Breiding,\:
John Cobb,\:
Aviva K. Englander,\:
Nayda Farnsworth,\\
Jonathan D. Hauenstein,\:
Oskar Henriksson,\:
David K. Johnson,\\
Jordy Lopez Garcia,\:
Deepak Mundayur
}

\date{}

\maketitle

\vspace{1em}

\begin{abstract}
\noindent
Many hypersurfaces in algebraic geometry, such as discriminants, arise as the projection of another variety. The real complement of such a hypersurface decomposes into connected components. In this paper, we propose a new method for computing these components. 
Existing methods require the explicit equation of the hypersurface as input.
However, computing this equation by elimination can be computationally demanding or even infeasible.
Our approach instead derives from univariate interpolation 
by computing the intersection of the hypersurface with a line.  
Such an intersection may be computed using so-called pseudo-witness sets \emph{without} computing a defining equation for the hypersurface. 
We implement our approach in a forthcoming Julia package and demonstrate, on several examples, that the resulting algorithm accurately recovers all components of the
real complement of the hypersurface.

\medskip


\medskip


\vspace{2em}

\end{abstract}

\section{Introduction}\label{sec:introduction}

Many hypersurfaces in algebraic geometry arise as the Zariski closure of the projection of another variety. 
This paper presents a new method for computing the connected components, or \emph{regions}, of the real complement $\bR^n\smallsetminus\sH$ of a hypersurface $\sH\subset\bC^n$ defined over $\bR$ without requiring an explicit defining
equation for $\sH$. 
By \textit{computing} the regions, we mean finding a representation of them that allows one to obtain sample points from each region, determine the number of regions, and perform membership tests.

Several algorithms have been developed to compute the regions assuming a known defining polynomial $h$ for~$\sH $; see for example~\cite[Section 15]{BasuPollackRoyBook} and~\cite{CannyRoadmaps,SafeySchost03}. One such approach~\cite{HypersurfaceRegions, SmoothConnectivity,HongConnectivity,HRSS} uses so-called \emph{routing functions} of the form
\begin{equation}\label{eq:Routing}
r(x) = \frac{\vert h(x)\vert }{(1+\lVert x-c\rVert^{2})^{e}} = \frac{\vert h(x)\vert }{(1+(x_{1}-c_{1})^{2}+\cdots+(x_{n}-c_{n})^{2})^{e}},
\end{equation}
where $c\in \bR^{n}$ is generic and $e$ is a positive integer with $2e > \deg h$. This ensures that $r$ vanishes on $\sH$ and at infinity, so that $r$ has at least one critical point in every region of $\bR^n\smallsetminus\sH$. Computing these critical points and connecting them via gradient flow of $r$ yields a graph whose connected components represent the regions of $\bR^n\smallsetminus\sH$. This approach has found applications in, for instance, convex geometry \cite{BrandenburgMeroni2025} or ecology \cite{CummingsDahlinGrossHauenstein2025,CelikHaasScholtenWangZucal2025}, and has been implemented in the Julia package \texttt{HypersurfaceRegions.jl} \cite{HypersurfaceRegions}.

Many hypersurfaces of interest in applications arise as the projection of some variety $X=V(F)$ with known defining equations. In this case, obtaining a defining polynomial~$h$ of the hypersurface $\sH$ requires a computationally expensive symbolic elimination, which even for moderate problems is often out of reach with current technology, and thus prohibits the direct use of routing functions. 

One of the main observations of this work is that a \emph{pseudo-witness set} representation \cite{WitnessSetsOfProj} of $\sH$ allows us to carry out the following tasks \emph{without} explicitly computing the square-free defining polynomial $h$, which is unique up to scaling:
\begin{enumerate}[label=(\arabic*)]
\item compute the degree $\deg h$;
\item evaluate $\log \lvert h(x)\rvert$ at points $x\in\mathbb R^n\smallsetminus\sH$ up to an additive constant;
\item evaluate the gradient $\nabla \log \lvert h(x)\rvert$ at points $x\in\mathbb R^n \smallsetminus \sH$; and
\item evaluate the Hessian $\Hess(\log \lvert h(x)\rvert)$ at points $x\in\mathbb R^n \smallsetminus \sH$.
\end{enumerate}

These operations are enough to define a routing function $r$,
use homotopy continuation \cite{bates2024numericalnonlinearalgebra, BertiniBook,SommeseWamplerBook} to determine critical points of $r$, and run gradient flow. The key idea behind pseudo-witness sets is to replace the explicit defining polynomial $h$ with numerical data from the intersection of $\sH$ with a generic line, computed by lifting the intersection problem to $X$.

Pseudo-witness sets have previously been used for membership tests \cite{HauensteinSommese2013}, finding generators of certain elimination ideals \cite{HauensteinIkenmeyerLandsberg2013}, and, more generally, numerical implicitization \cite{ChenKileel2019}. Recently, in~\cite{ImplicitPolynomials}, it was proposed that pseudo-witness sets can be used for evaluating and computing directional derivatives of the defining polynomial of a hypersurface that arises through projection. Here, we build on this approach to show that pseudo-witness sets can be used to compute gradients and  Hessians (\Cref{thm:gradient_and_hessian}). This gives rise to an algorithm (\Cref{alg:complete}) for computing the regions of $\bR^n\smallsetminus\sH$ that circumvents the need for symbolic elimination.

\begin{figure}[htbp]
    \centering
    \includegraphics[width=0.55\textwidth]{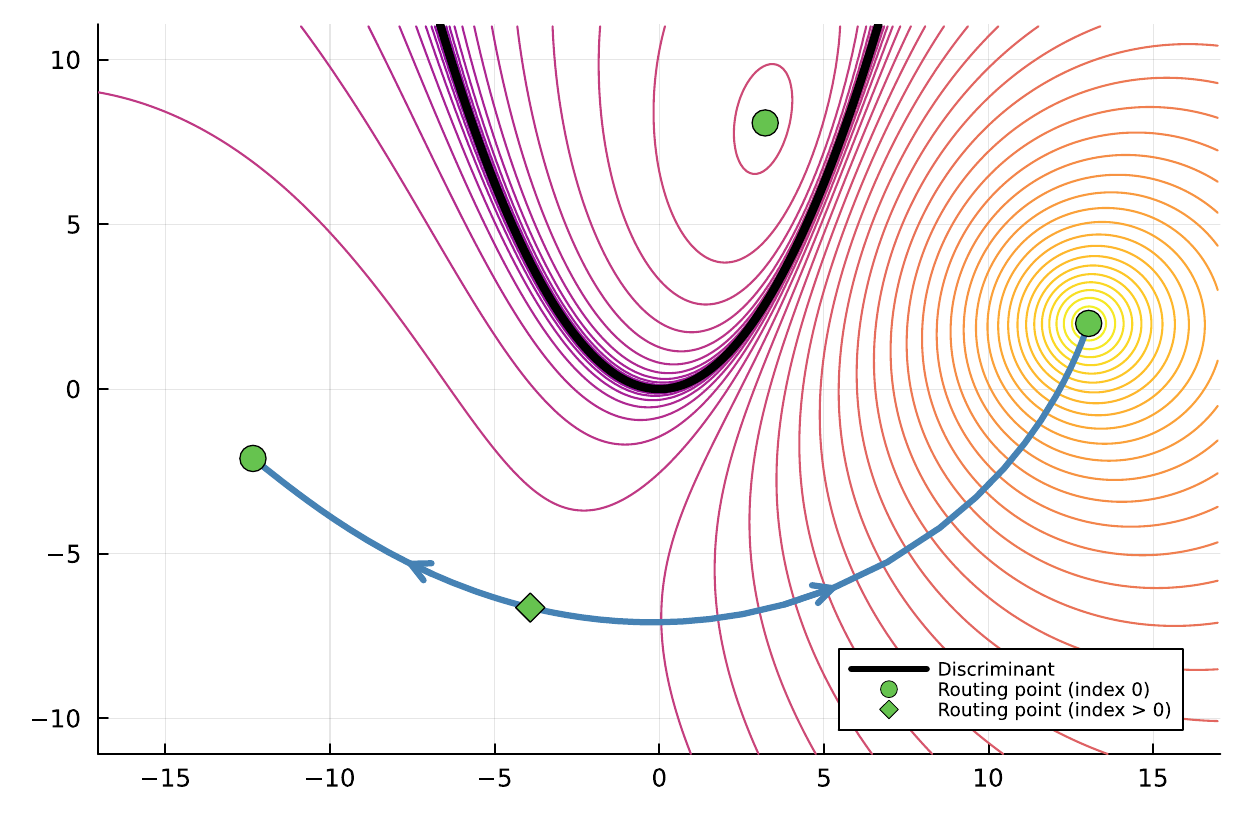}
    \caption{The real part of the quadratic discriminant $V(h)$ is the black curve. There are two regions: one above and one below the black curve. Our algorithm represents these regions by critical points (green) of a routing function (with level sets illustrated with colors between magenta and yellow). There is one critical point in the top region and three critical points in the bottom region. The latter three points are connected by gradient flow (the blue curves). The critical points and the flow trajectories were computed without direct access to the polynomial $h = a^2 -4b$.}
    \label{fig:quadratic-discriminant}
\end{figure}

\begin{example}\label{ex:quadratic_discriminant}
An important class of hypersurfaces that are obtained through projection is \emph{discriminants} of parametric polynomial systems (see, e.g., \cite{gkz1994} and \Cref{sec:examples}). As a running example, we will consider one of the simplest instances: the discriminant of the quadratic polynomial 
\begin{equation*}
f(a,b,z) = z^2 + az + b,  
\end{equation*}
where we view $(a,b)$ as parameters, and $z$ as the variable. We let $X\subset\bC^3$ be the variety of triples $(a,b,z)$ such that $z$ is a double root of $f(a,b,z)$, i.e., the zeros of the system 
\begin{equation}\label{eq:QuadDisc1}
F(a,b,z) = 
\begin{pmatrix}
f(a,b,z) \\[0.5ex] f'(a,b,z)
\end{pmatrix}
= 
\begin{pmatrix}
z^{2} + az + b \\[0.5ex] 2z + a
\end{pmatrix}.
\end{equation}
This yields an irreducible curve in $\bC^3$, and the closure of its projection to the parameter space $\mathbb C^2$ is the discriminant hypersurface
\begin{equation}\label{eq:QuadDisc}
\sH = \bigl\{(a,b)\in\bC^2\suchthat h(a,b):=a^2 - 4b=0\bigr\}.
\end{equation}
The real complement $\bR^{2}\smallsetminus \sH$ of the hypersurface consists of the two regions 
\begin{equation*}
\{(a,b) \in \bR^{2}  \mid h(a,b) >0\} \quad \text{and} \quad \{(a,b) \in \bR^{2} \mid h(a,b)<0\}. 
\end{equation*}
Our algorithm first uses a pseudo-witness set to determine that the degree of $h$ is $2$, picks a random point $c = (c_{1},c_{2})\in\bR^2$, and constructs the routing function
$$r(a,b) = \frac{\vert h(a,b)\vert }{(1+ (a-c_1)^2+(b-c_2)^2)^2}.$$ 
Using $c=(13,2)$, our algorithm finds four real critical points, and using gradient flow, it determines that three of them belong to one region, whereas the fourth belongs to another region; see \Cref{fig:quadratic-discriminant}. During this computation, our algorithm never has direct access to the polynomial $h$. 
\end{example}

The remainder of this paper is organized as follows. \Cref{sec:background} reviews background on routing functions, gradient roadmaps, and homotopy continuation for the case when the defining polynomial~$h$ of $\mathcal{H}$ is known. \Cref{sec:unknown_h} describes how to extend this approach to the case when $h$ is unknown, through the concept of pseudo-witness sets, and ends with a high-level description of our method (\Cref{alg:complete}). 
In \Cref{sec:examples}, we illustrate how our method can be applied to discriminant varieties to explore how the real or nonnegative root count of a parametrized polynomial system varies throughout the parameter space.
The paper concludes in~\Cref{sec:Conclusion}.

\section{The case of a known defining polynomial}
\label{sec:background}

This section collects background information on computing the regions of $\bR^{n}\smallsetminus\sH$ for a hypersurface $\mathcal{H}=V(h) \subset \bC^{n}$ with known defining polynomial $h$. In particular, we recall the concept of routing functions (\Cref{subsec:Routing}) 
and gradient roadmaps (\Cref{subsec:RoutingConnectivity}), as well as homotopy continuation and monodromy solving (\Cref{sec:HCMonodromy}). 

\subsection{Routing points and routing functions on \texorpdfstring{$\bR^{n}$}{Rn}}\label{subsec:Routing}

This section is based on \cite[Sections 2 and 3]{SmoothConnectivity}. Let $h\in\bR[x_1,\ldots,x_n]$ be a polynomial defining a hypersurface $\sH = V(h)\subset \bC^{n}$. The real complement $\bR^{n}\smallsetminus \sH$ consists of connected components, which we call \emph{regions}, on which the sign of $h$ is constant. To compute these regions, we use the routing function
\begin{equation}\label{eq:routing-function}
    r(x) = \frac{\vert h(x)\vert }{q(x)^{e}},\quad q(x)  = 1 + \sum_{i=1}^{n}(x_{i}-c_{i})^{2},
\end{equation}
where $c = (c_{1},\dots,c_{n}) \in \bR^{n}$ is a generic point and $e$ is a positive integer satisfying $2e > \deg h$. As discussed in \Cref{sec:introduction}, the idea is to compute the critical points of $r$. 

\begin{definition}\label{def:routing_points}
The \emph{routing points} of the routing function $r$ from \eqref{eq:routing-function} are given by the set 
\begin{equation*}
    Z(r) \vcentcolon = \bigl\{y\in\bR^{n}\mid r(y)\neq 0,\ \nabla r(y)=0\bigr\}.
\end{equation*}
A routing point $y\in Z(r)$ is called {\em nondegenerate} if the Hessian matrix $\Hess(r)(y)$ is invertible. The {\em index} of a nondegenerate routing point $y$ is the number of \emph{positive} eigenvalues of $\Hess(r)(y)$ counted with multiplicity. Nonzero eigenvectors of $\Hess(r)(y)$ corresponding to positive eigenvalues are called {\em unstable eigenvector directions}. 
\end{definition}

In our convention, the index of a nondegenerate routing point $y$ counts the positive eigenvalues of $\Hess(r)(y)$ rather than the negative eigenvalues counted by the usual Morse index. 
The next theorem summarizes the properties of the routing function $r$ from \cite{SmoothConnectivity} that allow us to use it for computing the connected components of $\bR^{n}\smallsetminus \sH$.

\begin{theorem}[\hspace{-0.04em}{\cite[Proposition~3.2, Theorem~3.4, Proof of Theorem~3.8]{SmoothConnectivity}}]
\label{thm:generic_properties_of_routing_functions}
Let $r$ be as in \eqref{eq:routing-function}, with  $e>\deg(h)/2$. Then, there exists a nonempty Zariski open subset $\mathcal{U}\subset\mathbb{R}^n$ such that for every $c=(c_1,\ldots,c_n)\in\mathcal{U}$, the following properties hold.
\begin{enumerate}[itemsep=0pt,parsep=0pt]
    \item The restriction of $r$ to the complement $\bR^n\smallsetminus\sH$ is a Morse function.
    \item For all $\epsilon >0$, there exists $A > 0$ such that if $x\in \bR^{n}$ with $\lVert x\rVert\ge A$, then $\lvert r(x)\rvert < \epsilon$.
    \item The set $Z(r)$ of routing points of $r$ 
    is finite.
    \item Each routing point is nondegenerate.
    \item For each $\alpha >0$,  there is at most one $y\in Z(r)$ with $r(y)=\alpha$.
    \item The norms of $r$, $\nabla r$, and $\Hess(r)$ are bounded on $\bR^{n}\smallsetminus \sH$.
    \item Each region of $\bR^{n}\smallsetminus \sH$ contains at least one routing point of $r$ of index~$0$.
\end{enumerate}
\end{theorem}

\begin{example}\label{Ex: quadratic discriminant4}
We continue our running example. Recall from \Cref{ex:quadratic_discriminant} the discriminant polynomial 
$h(a,b) = a^{2}-4b$.  
Taking $e=2$ and $(c_{1},c_{2}) = (13,2)$, we obtain 
\begin{equation*}
r(a,b) = \frac{\vert a^{2}-4b\vert }{\big(1+(a - 13)^{2}+(b -2)^{2}\big)^{2}},    
\end{equation*}
which vanishes at infinity and satisfies the conditions of \Cref{thm:generic_properties_of_routing_functions}.  
Thus $r$ is a routing function on $\bR^{2}$ and it has four (nondegenerate) routing points approximately at
\begin{equation*}
(-12.339, -2.107),\quad (-3.918, -6.636),\quad (13.040, 1.994),\quad \text{and} \quad (3.217, 8.083) 
\end{equation*}
with indices 0, 1, 0, and 0, respectively. Each region of  $\bR^{2}\smallsetminus \sH$ contains a routing point of index zero as guaranteed by \Cref{thm:generic_properties_of_routing_functions}; see \Cref{fig:quadratic-discriminant}.
\end{example}

Although the routing function $r$ is defined using the absolute value  of $h$, its logarithmic gradient admits a rational extension to complex points. Observe that the sign of $h$ is constant on each region of $\bR^{n}\smallsetminus \sH$, so 
\begin{equation*}
\nabla \log \lvert h(x)\rvert = \frac{\nabla h(x)}{h(x)}.     
\end{equation*}
Thus, for $x\in \bR^{n}\smallsetminus \sH$, the logarithmic gradient
\begin{equation*}
    \nabla \log r(x) = \nabla \log \frac{\lvert h(x)\rvert}{q(x)^{e}}
\end{equation*}
is the restriction to the real points of the rational map 
\begin{equation}\label{eq:Phi}
\Phi(x) := \frac{\nabla h(x)}{h(x)} - \frac{2e(x-c)}{q(x)},    
\end{equation}
defined on $\bC^{n}\smallsetminus(\sH \cup V(q))$. 

Since $r(x) > 0$ for $x\in \bR^{n}\smallsetminus \sH$, we have 
\begin{equation*}
\nabla r(x) = r(x)\Phi(x). 
\end{equation*}
Thus, the real zeros of $\Phi$ are the routing points. Moreover, at a routing point $y$,
\begin{equation*}
    J_{y}\Phi(y) = \Hess(\log r)(y) = \frac{\Hess(r)(y)}{r(y)}.
\end{equation*}
Since $r(y) > 0$, the eigenvalues of $J_{y}\Phi(y)$ are positive scalar multiples of those of $\Hess(r)(y)$, with the same eigenvectors. Thus, nondegeneracy, index, and unstable eigenvector directions are unchanged by passing from $r$ to $\log r$. Moreover, 
\begin{equation*}
    \nabla \log r(x) = \frac{\nabla r(x)}{r(x)}, 
\end{equation*}
so the gradient flows of $r$ and $\log r$ have the same trajectories, up to reparameterization.

\subsection{Gradient flow and connectivity via routing functions}\label{subsec:RoutingConnectivity}
A routing function $r$ together with its routing points $Z(r)$ from \Cref{subsec:Routing} constitutes a powerful numerical representation of the complement $\bR^{n}\smallsetminus \sH$. The properties of $r$ allow us to determine which of the routing points belong to the same region. This, in turn, yields a membership test for each of the regions, as well as the Euler characteristic of each region.

The key to connecting the routing points that belong to the same region is the \emph{Mountain Pass Theorem} \cite{mountain-pass-1,mountain-pass-2,mountain-pass-3}, which (in our setting) implies that all routing points of index zero belonging to the same region can be connected via gradient flow emanating from routing points of index one; see~\cite{mountain-pass-2},~\cite[Thm. 3]{mountain-pass-3} and the proofs of~\cite[Proposition 4.2, Thm. 4.4]{SmoothConnectivity}. This gives rise to \Cref{alg:connected-components}, which is an adaptation of \cite[Algorithm~1]{SmoothConnectivity}. The algorithm is illustrated by \Cref{fig:quadratic-discriminant}, where the three routing points that belong to the same region are connected via a single routing point of index one. 

\medskip
\begin{algorithm}[htbp]
\caption{Connected components}
\label{alg:connected-components}
\SetAlgoLined
\KwIn{A routing function $r$ on $\bR^{n}$ with numerator $h$ satisfying the assumptions of \Cref{thm:generic_properties_of_routing_functions}, and its routing points $Z(r)=\{y_1,\ldots,y_m\}$.}
\KwOut{A partition $\{C_{1},\dots,C_{s}\}$ of the set of routing points $Z(r)$ of $r$ on $\bR^{n}$ corresponding to the connected components of $\bR^{n}\smallsetminus V(h)$.} 
\BlankLine

\BlankLine
\textbf{(1) Initialize a graph.}
Construct a graph $G$ with vertex set $Z(r)$ and no edges.

\BlankLine
\textbf{(2) Add edges via unstable directions.} 

\BlankLine

\For{$j=1,\dots, m$}{
    \If{the index of $y_{j}$ is positive}{Choose an unstable eigenvector $v$ for $\Hess(r)(y_{j})$.\\
    Compute the limit routing point $y_{w_{+}}$ from $y_{j}$ in the direction $v$ by gradient flow with respect to $r$ (equivalently, $\log r$).\\
    Add the edge $\{y_j,y_{w_{+}}\}$ to $G$. \\
    Compute the limit routing point $y_{w_{-}}$ from $y_{j}$ in the direction $-v$ by gradient flow with respect to $r$ (equivalently, $\log r$).\\
    Add the edge $\{y_j,y_{w_{-}}\}$ to $G$.}
    }
\BlankLine
\textbf{(3) Transitive closure.}
Let $C_{1},\dots,C_{s}$ be the connected components of $G$.

\BlankLine
\textbf{(4) Return} $\{C_{1},\dots,C_{s}\}$.
\end{algorithm}

Once we have partitioned the routing points into subsets $C_1,\ldots,C_s\subset Z(r)$ as in \Cref{alg:connected-components}, we can use each $C_i$ as a representation of the region of $\bR^{n}\smallsetminus \sH$ that its routing points belong to. For instance, it can be used for membership testing: any point $z\in$  $\bR^{n}\smallsetminus \sH$ that is not a routing point belongs to the region represented by $C_i$ if and only if gradient flow with respect to $r$ from $z$ gives a path that converges to a routing point in $C_i$  \cite[Proposition~4.1]{SmoothConnectivity}. Furthermore, the Euler characteristic of the region represented by $C_i$ is given by $\sum_{j=0}^n (-1)^j\mu_j$, where $\mu_j$
is the number of routing points of index $j$ in $C_i$ \cite[Theorem~3.8]{SmoothConnectivity}.

\subsection{Homotopy continuation and monodromy solving}\label{sec:HCMonodromy}

In the previous two sections we explained how to use a routing function $r$ for computing regions of $\bR^{n}\smallsetminus \sH$. The key is computing routing points. For numerical reasons, we work with $\log r$ instead of $r$; as discussed at the end of~\Cref{subsec:Routing}, this does not change the routing points, their indices or unstable eigenvector directions, or the gradient flow trajectories. Moreover, to compute the routing points, we first solve for the complex zeros of $\Phi$ defined in~\eqref{eq:Phi} using \textit{homotopy continuation}.

Homotopy continuation is a foundational numerical method in algebraic geometry. For an introduction, we refer the reader to \cite{BertiniBook,SommeseWamplerBook}. All homotopies considered here are analytic functions
\begin{equation*}
\begin{split}
    H \colon \mathcal{U} \times \bC & \longrightarrow \bC^{n}\\
    (x,t) &\longmapsto H(x,t),
\end{split}
\end{equation*}
for an open domain $\mathcal{U}\subseteq\mathbb{C}^n$. Given a {\em start point} $x^{*}\in \bC^{n}$ satisfying
\begin{equation}\label{davidenko_condition}
    H(x^{*},1) = 0,\quad \det J_{x}H(x^{*},1)\neq 0,
\end{equation}
where $J_{x}H$ denotes the Jacobian of $H$ with respect to $x$, the Implicit Function Theorem guarantees a smooth local path $x(t)$ near $t=1$ such that $H(x(t),t)\equiv 0$ and $x(1) = x^{*}$. Following~\cite{Regeneration}, this path is {\em trackable} if all points along it are nonsingular solutions of $H(\bullet, t) = 0$ for $t\in (0,1]$.
A trackable path can be numerically followed using a predictor--corrector scheme to integrate the Davidenko differential equation
\begin{equation}\label{eq:Davidenko}
\dot{x}(t) = -J_{x}H(x,t)^{-1}\, J_{t}H(x,t), \quad H(x(t),t)=0,
\end{equation}
for $t \in (0,1]$. If the limit $\lim_{t\to 0}x(t)$ exists in $\mathcal{U}$, the path is said to {\em converge} to an {\em endpoint} $x(0)$ satisfying $H(x(0),0) = 0$; otherwise, the path {\em diverges}.

\begin{example}\label{Ex: quadratic discriminant5}
    Let $F(a,b,z)$ be as in \Cref{ex:quadratic_discriminant}. We wish to intersect the curve $V(F)\subseteq\bC^3$ with the plane $V(b-1)$ by solving $F(x_1,1,x_2)=0$. For this, we consider the homotopy
    \begin{equation*}
        H(x,t) = (1-t)F(x_{1},1,x_{2}) + \gamma \,t\begin{pmatrix}x_{1}^{2} - 1 \\[2pt] x_{2} - 1\end{pmatrix}= 0,
    \end{equation*}
    where $\gamma$ is a generic complex number so that  $J_{x}H(x,t)$ is nonsingular along $t\in(0,1]$. Explicitly,
    \begin{equation*}
        H(x,t) = (1-t) \begin{pmatrix}x_{2}^{2} + x_{1}x_{2} + 1 \\[2pt] 2x_{2} + x_{1}\end{pmatrix} + \gamma\, t\begin{pmatrix}x_{1}^{2} - 1 \\[2pt] x_{2} - 1\end{pmatrix}.
    \end{equation*}

    The start points $(1,1)$ and $(-1,1)$ are nonsingular solutions of $H(x,1)=0$. Both corresponding paths are trackable and converge to the solutions $(2,-1)$ and $(-2,1)$ of $H(x,0)=0$. This shows that the intersection of the zero set of $F$ with the plane $V(b-1)$ is $\{(2,1,-1), (-2,1,1)\}$.
\end{example}
When only a subset of the solutions to a system of polynomials or rational functions is known, monodromy provides a method for finding additional solutions~\cite{TensorDecomp,MonodromyDecomp,Duff2018monodromy,MonodromySolver,MonodromyStats}. In this approach, we view our target system as a member of a parametric family of systems, and let the parameters vary continuously along a loop in parameter space, while numerically tracking a known solution. When the loop returns to its starting point, the endpoint of the tracked path is again a solution of the original system, but not necessarily the same one. This induces an action of the fundamental group of the complement of the branch locus on the fiber of solutions over a base parameter, called the {\em monodromy action}. The image of this action inside the permutation group of the solutions is the {\em monodromy group}. See~\cite{GaloisMonodromy} for further discussion of monodromy groups and their numerical computation. A transitive monodromy group implies that all (complex) solutions of the system can be reached from a single start solution via finitely many monodromy loops.

\begin{example}\label{ex:QuadDiscMonodromy}
In~\Cref{Ex: quadratic discriminant5}, we computed the intersection of the curve $V(F)$ with the plane $V(b-1)$ by solving $F(x_1,1,x_2)=0$. We can vary that plane in a monodromy loop:
for a loop $\Gamma\colon [0,1] \to \bC$ with~$\Gamma(0)=\Gamma(1)=1$, consider the homotopy
$
        H(x,t) = F(x_{1},\Gamma(t),x_{2}) = 0
$
    with start point $(2,-1)$. If $\Gamma(t) = 3/4 + (1/4)e^{2\pi i(1-t)}$, the corresponding path is trackable and returns to $(2,-1)$. In contrast, for the loop $\Gamma(t) = 1/4 + (3/4)e^{2\pi i(1-t)}$, the path ends at $(-2,1)$. Hence, the monodromy action exchanges these two solutions. The reason why the second loop gives a new solution, while the first does not, is that the second loop encircles a so-called \emph{branch point}.
\end{example}

In our case, we embed our target system $\Phi(x) = 0$ into the family $\Phi(x) -\eta = 0$ for parameters $\eta\in\bC^n$. As a start solution, we take a generic point $x_0\in\bC^n\smallsetminus (\sH \cup V(q))$, which by design is a solution of $\Phi(x) - \eta_{0}=0$ if we set $\eta_{0} = \Phi(x_{0})$. Since the \emph{incidence variety} 
\begin{equation*}
\{(x,\eta)\in (\bC^{n}\smallsetminus(\sH \cup V(q)))\times \bC^{n} \mid \Phi(x)-\eta=0\}
\end{equation*}
is a graph of a rational map over an irreducible Zariski open subset, and hence irreducible, it follows that the monodromy action is transitive \cite{Duff2018monodromy}. Consequently, we can find all solutions of $\Phi(x) - \eta_{0}=0$ by tracking monodromy loops until some (heuristic) stopping condition is reached. We then track the solutions along $H(x,t) = \Phi(x) - t\eta_0$  from $t=1$ to $t=0$; the finite real endpoints in $\bR^{n}\smallsetminus \sH$ are the routing points.

\subsection{Summary of the algorithm for known defining polynomial}
\label{subsec:summary_known_h}

We discussed all the ingredients needed to compute the regions of $\bR^{n}\smallsetminus \sH$ for a known polynomial $h$, and we summarize this as \Cref{alg:complete_known_h}, which is an adaptation of  \cite[Algorithm~1]{SmoothConnectivity}. This is also the foundational algorithm behind the Julia package \texttt{HypersurfaceRegions.jl} \cite{HypersurfaceRegions}. 
The goal of the rest of the paper is to adapt this to the case when the defining polynomial $h$ is unknown, and the hypersurface $\mathcal{H}$ instead arises as the projection of a known variety.

\begin{algorithm}[htbp]
\caption{Computing regions of a hypersurface with known defining polynomial}
\label{alg:complete_known_h}
\SetAlgoLined
\KwIn{A nonzero polynomial $h\in\bR[x_1,\ldots,x_n]$ defining a hypersurface $\sH$.}
\KwOut{The regions of the real complement $\bR^{n}\smallsetminus \sH$.} 
\BlankLine
\textbf{(1) Construct a routing function.}
Randomly select $c\in\bR^n$ and take $e>\deg(h) / 2$. This defines a routing function $r(x) = \vert h(x)\vert /q(x)^{e}$ as in \eqref{eq:routing-function}. 

\BlankLine
\textbf{(2) Compute routing points.}
Compute all solutions to $\Phi(x) = 0$ defined as in~\eqref{eq:Phi} through monodromy, as described in \Cref{sec:HCMonodromy}. 
Retain the real solutions. These are the routing points. 

\BlankLine
\textbf{(3) Compute regions.}
Run \Cref{alg:connected-components}.
\end{algorithm}

\section{The case of an unknown defining polynomial}
\label{sec:unknown_h}

In this section, we treat the main problem of interest in the paper, namely, to compute the regions of the complement $\bR^k\smallsetminus \mathcal{H}$ of a hypersurface  $\mathcal{H}\subset\bC^k$ with an unknown defining polynomial $h\in\mathbb{R}[p_1,\ldots,p_k]$, defined as the closure of the projection of some known variety $X \subset \bC^{n}$. More precisely, we consider an irreducible component 
$X \subset V(F)\subset\bC^{n}$ for a known system 
$F\in (\bR[p_1,\ldots,p_k,z_1,\ldots,z_{n-k}])^{m}$
such that
\begin{equation*}
\mathcal{H}=\overline{\pi(X)},
\end{equation*}
where $\pi\colon\bC^k\times\bC^{n-k}\to\bC^k$ is the projection $\pi(p,z)=p$, and the overline denotes Zariski closure. 
In what follows, we will let $h$ be a polynomial of minimal total degree such that $\mathcal{H}=V(h)$ (this defines $h$ uniquely up to scaling).

The key observation is that in order to carry out  \Cref{alg:complete_known_h}, we do not need to have an explicit expression for $h$. Instead, we just need the following.
\begin{itemize}
    \item We need to know $\deg h$ in order to choose the denominator of the routing function in step~(1).
    \item We need to be able to evaluate $(\nabla h)/h$ and $J_{p}(\nabla h/h)$ in order to evaluate the rational map $\Phi$ and its Jacobian $J_{p}\Phi$ and run homotopy continuation in step~(2) and gradient flow in step~(3).
\end{itemize}
It turns out that both these pieces of information about $h$ are available through a  \emph{pseudo-witness set} of $\mathcal{H}$. We explain this concept and how it gives the degree in \Cref{sec:pseudo-witness}. We then proceed to derive algorithms for evaluating $(\nabla h)/h$ and $J_{p}((\nabla h)/h)$, and hence $\Phi$ and $J_{p}\Phi$, from the data of a pseudo-witness set in \Cref{subsec:evaluating}, which is the main theoretical contribution of the paper. Building on this, we give a strategy (\Cref{alg:complete}) for computing $\bR^k\smallsetminus \mathcal{H}$ in \Cref{subsec:summary}.

\begin{remark}\label{rmk:multiple_components}
For the case when $X$ consists of several irreducible components, we can represent  $\mathcal{H}$ by one pseudo-witness set for each component, each contributing one summand of $\Phi$. When several components of $X$ project to the same component of $\mathcal{H}$, it suffices to work with one of them.   
\end{remark}

\begin{remark}\label{extend_h}
It is often useful to consider the complement of
the union of $\sH$ with another hypersurface $V(g)$ 
where $g$ is a known polynomial, 
e.g., when imposing positivity conditions~\cite{SmoothConnectivity,CummingsDahlinGrossHauenstein2025}.  This situation can be handled by a \textit{logarithmic barrier}, i.e., by simply adding $\log(|g|)$ to $\log(r)$ above 
and adjusting $e$ to ensure that $2e>\deg(h)+\deg(g)$. The corresponding rational map is obtained by replacing $\Phi$ with $\Phi + (\nabla g)/g$ defined away from $\sH \cup V(g)\cup V(q)$.
\end{remark}

\subsection{Witness and pseudo-witness sets}\label{sec:pseudo-witness}
Witness and pseudo-witness sets are data structures for representing (positive-dimensional) algebraic varieties numerically. We give a brief introduction to this topic. 
General references for witness and pseudo-witness sets are the books~\cite{BertiniBook,SommeseWamplerBook}.

Let us first define the notion of a witness set. Let $F=(f_1,\ldots,f_m)\in(\bC[x_1,\ldots,x_n])^m$ be a system of $m$ polynomials and consider the corresponding variety $V(F)\subset\bC^n$. Suppose that $X\subset V(F)$ is a $d$-dimensional irreducible component of $V(F)$. A witness set represents $X$ by the intersection of $X$ with a general linear space $\mathcal{M}$ of codimension $d$; this intersection may be found by solving the system obtained by augmenting $F$ with $d$ equations of $\mathcal{M}$. For a sufficiently general choice of $\mathcal{M}$, the intersection $X\cap \mathcal{M}$ consists of $\deg X$ points. Below is the full definition. 

\begin{definition}\label{def:witness_set}
\begin{enumerate}\item 
Let $X \subset \mathbb{C}^n$ be an irreducible variety of dimension $d$.
A \textit{witness set} for $X$ is a triple $(F,\mathcal{M}, \mathcal{W})$, where $F$ is a system such that $X$ is an irreducible component of $V(F)$, $\mathcal{M}\subset \bC^{n}$ is a general linear space of codimension $d$ and $\mathcal{W} = X \cap \mathcal{M}$ is a finite set of points.  Here, $F$ is called a {\em witness system}, $\mathcal{M}$ is the {\em witness slice}, and $\mathcal{W}$ is the corresponding {\em witness point set}.
\item We call the witness set $(F,\mathcal{M}, \mathcal{W})$ \emph{reduced}, if $\mathrm{rank}(J_{x}F(x)) = n-d$ for all points $x\in\mathcal W$. Here, $J_{x}F(x)$ denotes the Jacobian of $F$ at $x$. 
\item If $X$ is of pure dimension $d$, but not necessarily irreducible, let $X_1,\ldots,X_\ell$ be its irreducible components. We call $(F, \mathcal M, \mathcal{W}_1\cup \cdots\cup \mathcal W_\ell)$
a witness set for $X$, where $\mathcal W_i = X_i\cap \mathcal M$.
\end{enumerate}
\end{definition}

The second item can be rephrased as follows: 
a witness set for an irreducible component $X$ of $V(F)$ is reduced when the primary ideal corresponding to $X$  in the primary decomposition of $\langle F\rangle \vcentcolon= \langle f_1,\ldots, f_m\rangle $ is radical. By taking $n-d$ generic linear combinations of $f_1,\ldots,f_m$, the points in $X\cap \mathcal{M}$ become nonsingular solutions of a square polynomial system, in such a way that homotopy continuation can be used to trace them to another witness set $X\cap\mathcal{M}'$ for a different slice $\mathcal{M}'$. In this way, reduced witness sets enable several numerical operations in homotopy continuation, including sampling points on a variety, testing membership, and computing intersections~\cite{bates2024numericalnonlinearalgebra,hauenstein2012whatisnag,IntersectingWitnessSets,sommese2004homotopies,sommese2004intrinsic,SommeseWamplerBook}. 
If a witness set $(F,\mathcal{M},\sW)$ 
is not reduced, one can utilize the deflation techniques of
\cite{IsosingularDeflation} at the cost of introducing new variables and equations.

\begin{figure}[htbp]
    \centering
    \begin{picture}(220,145)
    \put(0,0){\includegraphics[scale=0.09]{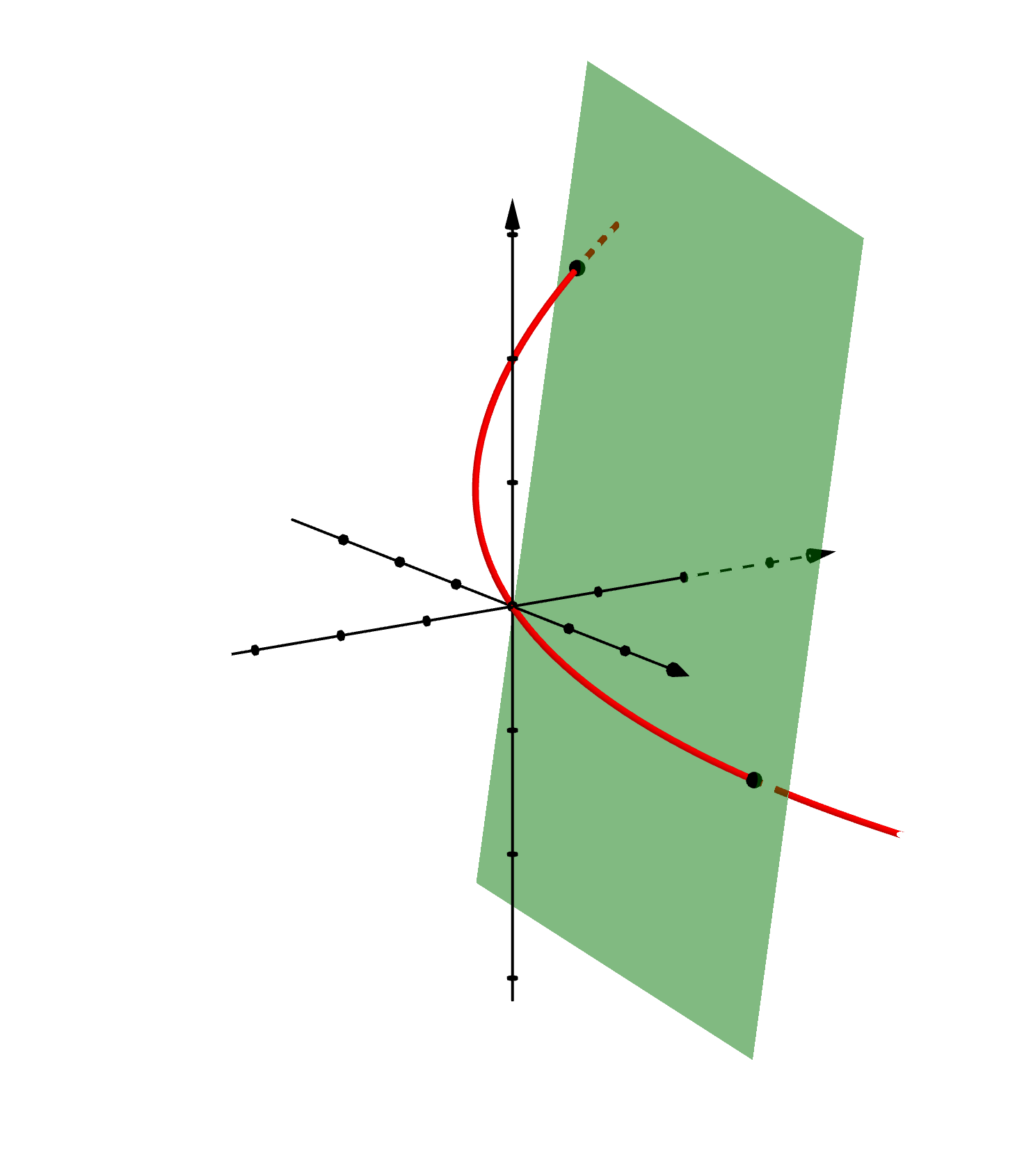}}
    \put(83,35){$\mathcal{M}$}
    \put(147,110){$\mathcal{W}$}
    \put(50,110){$X$}
    \put(145,113){\tikzmark{W-tail1}}
    \put(145,110){\tikzmark{W-tail2}}
    \put(81,117){\tikzmark{W1}}
    \put(100,55){\tikzmark{W2}}
    \end{picture}
    \begin{tikzpicture}[remember picture, overlay]
    \draw [line width=1.2pt, , ->, >=stealth] ({pic cs:W-tail1}) -- ({pic cs:W1});
    \draw [line width=1.5pt, ->, >=stealth] ({pic cs:W-tail2}) -- ({pic cs:W2});
    \end{tikzpicture}
    \caption{A witness set for the red curve $X = V(F)$ 
    with $F$ as in \Cref{ex:quadratic_discriminant2}. The witness set is given by the intersection of the red curve $X$ with the green linear space $\mathcal M$. The intersection consists of the two points labeled~$\mathcal W$.}
    \label{fig:WitnessQuad} 
\end{figure}

\begin{example}\label{ex:quadratic_discriminant2}
    We return to our running example from \Cref{ex:quadratic_discriminant}. 
    Let $F(a,b,z)$ be as in \eqref{eq:QuadDisc1} and recall that the zero set $V(F) \subset \bC^{3}$ of $F$ is an irreducible curve. Let $X = V(F)$ and intersect it with a general two-dimensional plane $\mathcal{M}$, say $V(a + 5b - z - 10)$. The intersection $\mathcal{M} \cap X$ is illustrated in \Cref{fig:WitnessQuad}. It consists of two points
    (in particular, this shows that $\deg X = 2$). The witness set for~$X$ is $(F,\mathcal M,\mathcal W)$ and it consists of~$F$, the plane $\mathcal M$ and the two intersection points in~$\mathcal{W} = X\cap \mathcal{M}$. The witness set $(F,\mathcal M,\mathcal W)$ is reduced.
\end{example}

\begin{example}
We give an example of a non-reduced witness set. Consider the case $n=2$ and the polynomial $F(x_1,x_2) = (x_1^2 + x_2^2 - 1)^2$. The zero set $X=V(F)$ is the circle. It has dimension~$d=1$. The intersection of $V(F)$ with a general line consists of $2$ points. 
Take, for instance, the line $\mathcal M=V(x_1)$. Then, $\mathcal W=  X\cap \mathcal M = \{(0,1), (0,-1)\}$. However, $J_{(x_1,x_2)}F(x_1,x_2) = 4(x_1^2 + x_2^2 - 1)(x_1,x_2)$, 
so that $J_{(x_1,x_2)}F(0,1) = (0,0)$ and $J_{(x_1,x_2)}F(0,-1) = (0,0)$ both have rank $0$, and not $n-d = 2-1 = 1$. The witness set $(F,\mathcal M,\mathcal W)$ is therefore \emph{not reduced}. This is caused by the fact that $x_1^2 + x_2^2 - 1$ appears with multiplicity 2 in $F$.
\end{example}

Next, we consider pseudo-witness sets. These represent the Zariski closure of images of projections of $X$. For this, let $F$ be in the variables $x=(x_{1},\dots, x_{n})$, with $p:=(x_{1},\dots, x_{k})$ denoting the first $k$ variables and $z:=(x_{k+1},\ldots,x_n)$ the last $n-k$ variables. We then consider the projection
\begin{equation}\label{projection}
    \pi \colon \bC^{n} \to \bC^{k}, \quad \pi(p,z) = p,
\end{equation}
and let 
\begin{equation}\label{def_H_formal}\sH = \overline{\pi(X)}\subset \bC^{k}\end{equation}
denote the Zariski closure of the image of $X$. 
When $\sH$ is a hypersurface whose defining polynomial is not explicitly known, we can still represent it numerically by a pseudo-witness set. The idea is to lift lines in $\mathbb C^k$ to a linear space in $\mathbb C^n$. This works as follows. Consider a general line $\mathcal{L}\subset \bC^{k}$ and a general linear space $\mathcal{L}' \subset \bC^{n-k}$ of dimension $n-d-1$. Then, the intersection of $\mathcal{L}\times \mathcal{L}'$ with~$X$ is a finite set $\mathcal{W}$. This set may be found by solving the system obtained by augmenting $F$ with $k-1$ equations of $\mathcal{L}$ and $d-k+1$ equations of $\mathcal{L}'$ and restricting to $X$. Projecting the points in $\mathcal{W} \subset \bC^{n}$ to the first $k$ coordinates yields the intersection $\pi(\mathcal W)=\sH\cap \mathcal L$. Note that a point in $\sH$ might have multiple preimages in $X$. 

This yields the notion of a \emph{pseudo-witness set}, which we summarize in the following definition. 
For simplicity, we only give the definition in the hypersurface case.

\begin{definition}\label{def:PWS}
\begin{enumerate}\item 
Let $X\subset \bC^n$ be an irreducible variety of dimension $d$, such that $\sH = \overline{\pi(X)}$ is a hypersurface.
A {\em pseudo-witness set}  for $\sH$ is a quadruple $(F,\pi,\mathcal{L}\times \mathcal{L}', \mathcal{W})$, where $F$ is a system such that $X$ is an irreducible component of $V(F)$, $\mathcal L$ is a general line and $\mathcal L'$ a general linear space  of dimension $n-d-1$, and  $\mathcal{W} = X\cap (\mathcal{L}\times \mathcal{L}')$.
\item We call the pseudo-witness set $(F,\pi, \mathcal{L}\times \mathcal{L}', \mathcal{W})$ \emph{reduced}, if $\mathrm{rank}(J_{x}F(x)) = n-d$ for all points $x\in\mathcal W$. Here, $J_{x}F(x)$ denotes the Jacobian of $F$ at $x$. 
\item 
If $X$ is of pure dimension $d$, but not necessarily irreducible, let $X_1,\ldots,X_\ell$ be its irreducible components
such that $\overline{\pi(X_j)}$ is a 
hypersurface.  
We call $(F, \pi, \mathcal{L}\times \mathcal{L}', \mathcal{W}_1\cup\cdots\cup \mathcal W_\ell)$
a pseudo-witness set for $\mathcal H = \overline{\pi(X)}$, where $\mathcal W_i = X_i\cap (\mathcal{L}\times \mathcal{L}')$.
\end{enumerate}
\end{definition}

As discussed above, once we want to use homotopy continuation and track a pseudo-witness set from one linear space $\mathcal{L}\times \mathcal{L}'$ to another $\hat{\mathcal{L}}\times \hat{\mathcal{L}}'$ we need to work with reduced pseudo-witness sets, so that the second condition in \eqref{davidenko_condition} is satisfied. 

\begin{example}\label{ex:quadratic_discriminant3}
We continue \Cref{ex:quadratic_discriminant,ex:quadratic_discriminant2}.
Recall that the discriminant $\sH$ is defined by the polynomial $h(a,b) = a^{2} - 4b$. It is obtained after eliminating $z$ from $F(a,b,z)$, where $F$ is as in~\eqref{eq:QuadDisc1}.
Consider the line 
$\mathcal L = \{ p + tv \mid t\in\mathbb C\},$
where $p=(0,2)$ and $v=(-2,3/5)$. We compute the intersection $\sH\cap \mathcal L$ using pseudo-witness sets for $X=V(F)$. Since $X$ is a curve ($d=1$) and since $n = 3$, we must take a general linear space $\mathcal L'$ of dimension $n-d-1=3-1-1=1$ in $\mathbb C$. There is only one option, namely $\mathcal L'=\mathbb C$. Computing the intersection of $X$ with $\mathcal L\times \mathcal L'$ yields 
\begin{equation*}
\mathcal{W}
= X\cap (\mathcal{L} \times \bC)
=
\Bigl\{
\Bigl(-\tfrac{3+\sqrt{209}}{5},\tfrac{109+3\sqrt{209}}{50}, \tfrac{3+\sqrt{209}}{10}\Bigr),
\Bigl(-\tfrac{3-\sqrt{209}}{5},\tfrac{109-3\sqrt{209}}{50}, \tfrac{3-\sqrt{209}}{10}\Bigr)
\Bigr\}.
\end{equation*}
Since $\pi(\mathcal{W}) = \sH \cap \mathcal{L}$, the intersection of the discriminant with $\mathcal{L}$ is given by taking the first two coordinates of these points. The pseudo-witness set $(F,\pi, \mathcal{L} \times \bC, \mathcal W)$ is reduced.
\end{example}

\begin{figure}[htbp]
    \centering
    \begin{subfigure}{0.45\textwidth}
        \centering
        \begin{picture}(200,190)
        \put(0,0){\includegraphics[scale=0.13]{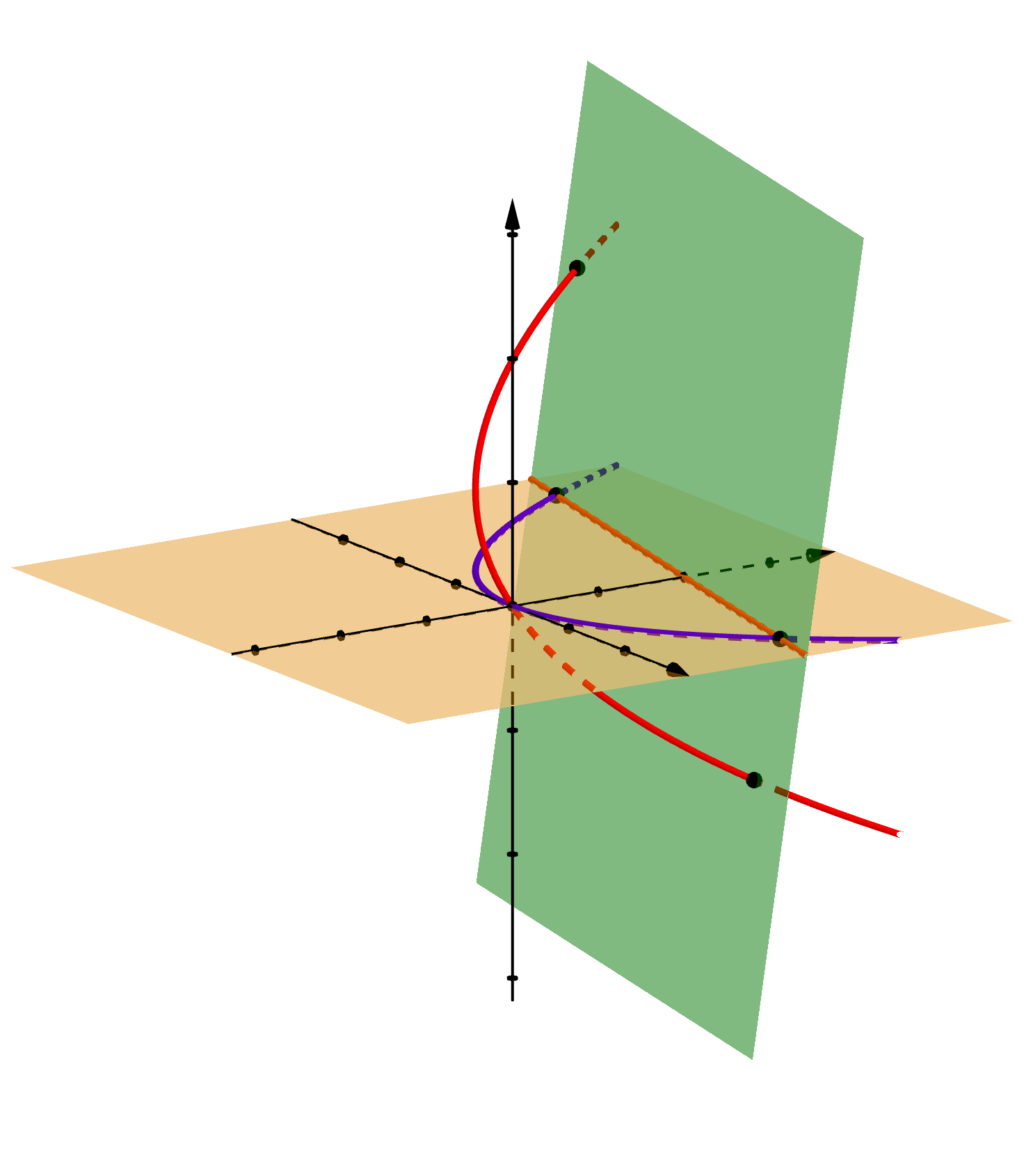}}
        \put(160,87){$\sH$}
        \put(120,120){$\mathcal{L}$}
        \put(128,150){$\pi(\mathcal{W})$}
        \put(55,135){$V(F)$}
        \put(135,145){\tikzmark{Wb-tail1}}
        \put(143,144){\tikzmark{Wb-tail2}}
        \put(108,129){\tikzmark{W1b}}
        \put(146,105){\tikzmark{W2b}}
        \end{picture}
        \begin{tikzpicture}[remember picture, overlay]
        \draw [line width=1.2pt, , ->, >=stealth] ({pic cs:Wb-tail1}) -- ({pic cs:W1b});
        \draw [line width=1.2pt, ->, >=stealth] ({pic cs:Wb-tail2}) -- ({pic cs:W2b});
        \end{tikzpicture}
        \vspace{-2em}
    \end{subfigure}
    \begin{subfigure}{0.45\textwidth}
        \centering
        \begin{picture}(150,190)
        \put(-20,0){\includegraphics[height=2in]{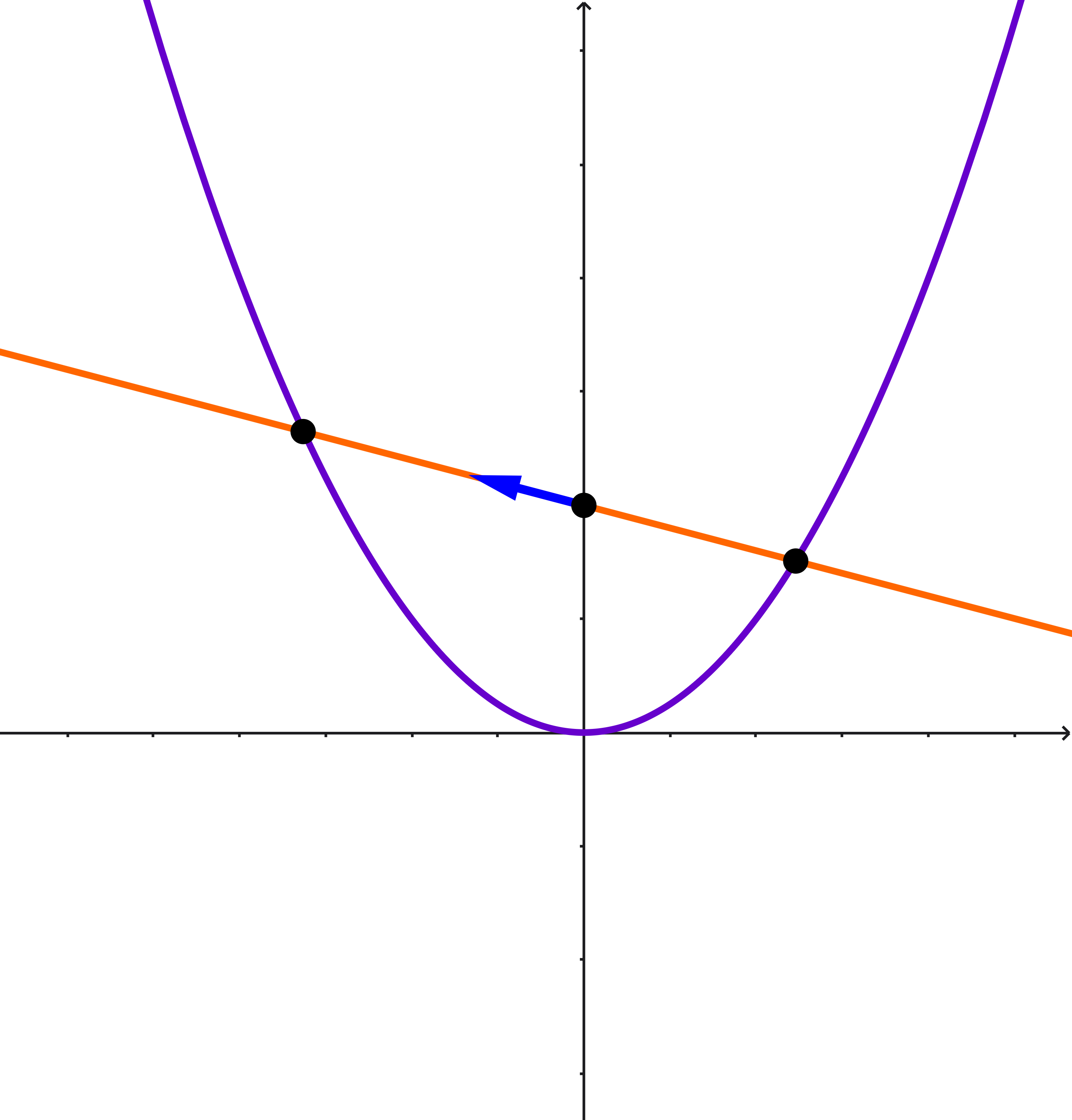}}
        \put(-15,85){$\mathcal{L}$}
        \put(108,110){$\sH$}
        \put(58,110){$\pi(\mathcal{W})$}
        \put(65,105){\tikzmark{piWb-tail1}}
        \put(75,104){\tikzmark{piWb-tail2}}
        \put(22,90){\tikzmark{piW1b}}
        \put(81,75){\tikzmark{piW2b}}
        \put(60,83){$p$}
        \put(40,70){$v$}
        \put(60,60){\tikzmark{btail}}
        \put(40,80){\tikzmark{bhead}}
        \end{picture}
        \begin{tikzpicture}[remember picture, overlay]
        \draw [line width=1.2pt, , ->, >=stealth] ({pic cs:piWb-tail1}) -- ({pic cs:piW1b});
        \draw [line width=1.2pt, ->, >=stealth] ({pic cs:piWb-tail2}) -- ({pic cs:piW2b});
        \end{tikzpicture}
    \end{subfigure}

    \caption{The left picture illustrates a pseudo-witness set for the purple projection $\sH$ of the red curve $V(F)$. While in \Cref{fig:WitnessQuad}, the green linear space $\mathcal M$ was general, here it must be a product space: $\mathcal M = \mathcal L \times \mathbb C$, where $\mathcal L$ is a general line; see \Cref{ex:quadratic_discriminant3}. The pseudo-witness set is given by the intersection $\mathcal{W}=V(F)\cap (\mathcal{L}\times\bC)$, and projection yields $\sH \cap \mathcal L$. The right picture shows the situation within the yellow plane. Here, $v \in \bR^{2}\smallsetminus \{(0,0)\}$ and $p \in \bR^{2}\smallsetminus \sH$.}
    \label{fig:PWSQuad}
\end{figure}

\subsection{Evaluating the gradient and Hessian}
\label{subsec:evaluating}

We now turn to the problem of evaluating $(\nabla h)/h$ and $J_{p}((\nabla h)/h)$ of the unknown defining polynomial~$h$ of $\sH \subset \bC^{k}$ at points $p\in \bC^{k}\smallsetminus \sH$ using the data of a pseudo-witness set $(F, \pi, \mathcal{L}\times\mathcal L', \mathcal{W})$ where $\mathcal{L}\subset\bC^k$ is a line. 

It is described in~\cite{ImplicitPolynomials}
how to utilize a pseudo-witness set to
evaluate and differentiate~$h$ along the line $\mathcal{L}$. This implies that $\nabla h(p)$ can be computed from $k$ linearly independent lines $\mathcal{L}$ through~$p$. Likewise, $\Hess(h)(p)$ can be computed through $k(k+1)/2$ lines through $p$. In this section, we prove that, in fact, a \emph{single} general line $\mathcal{L}$ through $p$ suffices to evaluate both $(\nabla h)/h$ and $J_{p}((\nabla h)/h)$. This drastically improves the complexity of the computation.

Let $\mathcal L\subset \bC^k$ be a general line that passes through a given point $p\in\bC^k \smallsetminus \mathcal{H}$ in the direction $b\in \bC^{k}$, and consider a general linear space $\mathcal L'$ of dimension $n-d-1$. Recall from \Cref{sec:pseudo-witness} that under the assumption of sufficient genericity, the line $\mathcal{L}$ meets $\sH$ transversely in $\deg h$ points and no intersection occurs at infinity. Let 
\begin{equation*}
\sH \cap \mathcal L = \{p+t_j b\mid j=1,\ldots,\deg h\}.
\end{equation*}
Then, we have that
\begin{equation}\label{eq:factor0}
    h(p+tb) = C(b)\cdot (t-t_1)\cdots (t - t_{\deg h})
\end{equation}
for some constant $C(b)$ that depends
on $b$.  
In fact, for generic $b$, 
$C(b)\neq0$ and $C(b)\cdot t^{\deg h}$ is the leading (highest-total degree) part of $h(p+tb)$. 
When $p$ and $b$ are real, setting $t = 0$ in~\eqref{eq:factor0}, taking absolute values and using logarithms gives
\begin{equation}\label{eq:factor0_log}
    \log |h(p)| = \log |C(b)| +\sum_{j=1}^{\deg h} \log |t_j| .
\end{equation}
Thus, by keeping $b$ constant, we can evaluate $\log |h(p)| $ up to an unknown constant $\log |C(b)| $.

A key component in what follows next is the gradient of $t_j$ with respect to $p$ and $b$. For computing these gradients, we choose parameterizations\enlargethispage{\baselineskip}
$$\mathcal{L} = \{p +t b \mid t\in \bC\}\quad\text{and}\quad \mathcal{L}' = \{v + Ay \mid y\in \bC^{n-d-1}\},$$
where $p \in \bC^{k}\smallsetminus \sH$ and $b\in \bC^{k}\smallsetminus\{0\}$, and $A\in\mathbb C^{(n-k)\times (n-d-1)}$ has full rank and $v\in \bC^{n-k}$. 

Since $p\notin\sH$, it will be helpful to work with the inverse $s \coloneq 1/t$ instead of $t$. A change of variables from $t$ to~$s$ yields the Laurent polynomial system 
$$G(s,y) = F(p + (1/s)b, v+Ay).$$
For each $j=1,\ldots,\deg h$, let $s_j=1/t_j$ and pick $y_j\in\bC^{n-d-1}$ such that $G(s_{j},y_{j})=0$. If the pseudo-witness set is reduced, $J_{(s,y)}G(s_j,y_j)$ has rank $n-d$, and it follows from the Implicit Function Theorem that locally around $(s_{j},y_{j})$ there exist smooth functions
\begin{equation}
\label{eq:definition_of_sj_and_yj}
s_j(p, b)\quad \text{and}\quad y_j(p, b),\qquad j=1,\ldots,\deg h,
\end{equation}
such that $G(s_j(p,b), y_j(p,b))=0$. 
Differentiating that equation with respect to $p$ and $b$ yields
\begin{equation}\label{eq:gradient-wrt-p-b}
\begin{aligned}
    J_{(s,y)}G
    \begin{pmatrix}
        \nabla_{p}s \\ J_{p}y
    \end{pmatrix} 
    + J_{p}G &= 0 \quad \text{and} \quad 
    J_{(s,y)}G
    \begin{pmatrix}
        \nabla_{b}s \\ J_{b}y
    \end{pmatrix} 
    + J_{b}G &= 0.
    \end{aligned}
\end{equation}
Differentiating the first system in~\eqref{eq:gradient-wrt-p-b} again with respect to each coordinate $b_{i}$ yields
\begin{equation}\label{eq:hessian-system}
J_{(s,y)}G
\begin{pmatrix}
    \partial_{b_{i}}(\nabla_{p}s) \\ \partial_{b_{i}}(J_{p}y)
\end{pmatrix}
+
R^{(i)}(p,b) = 0,
\end{equation}
where $R^{(i)}(p,b)$ depends only on $\nabla_{p}s, J_{p}y, \nabla_{b}s, J_{b}y$ and derivatives of the Jacobian matrices. Solving~\eqref{eq:hessian-system} yields $J_{b}(\nabla_{p}s)$. 

Differentiating~\eqref{eq:factor0} with respect to $t$ and dividing by $h(p+tb)$ gives
\begin{equation}\label{eq:directional-rational-h-complex}
    \frac{\frac{\mathrm{d}}{\mathrm{d}t}h(p+tb)}{h(p+tb)} = \sum_{j=1}^{\deg h}\frac{1}{t-t_{j}},
\end{equation}
which is valid over $\bC$ whenever the denominators are nonzero. This expression will be used to obtain formulas for the gradient and Hessian of $\log |h|$. 

\begin{proposition}
\label{prop:gradient}
For $p\in \bC^{k}\smallsetminus\sH$ and a generic direction $b \in \bC^{k}$,
\begin{equation*}
    \frac{\nabla h(p)}{h(p)} \cdot b = -\sum_{j=1}^{\deg h}s_{j}(p,b).
\end{equation*}
When $p$ and $b$ are real, the left-hand side agrees with the directional derivative of $\log \lvert h\rvert$ at $p$ in direction $b$.
\end{proposition}

\begin{proof}
Since 
\begin{equation*}
\frac{\mathrm{d}}{\mathrm{d}t} h(p +tb)\vert_{t=0} = \nabla h(p)\cdot b,
\end{equation*}
evaluating~\eqref{eq:directional-rational-h-complex} at $t=0$ gives
\begin{equation*}
\frac{\nabla h(p)}{h(p)}\cdot b = -\sum_{j=1}^{\deg h}\frac{1}{t_{j}} = -\sum_{j=1}^{\deg h}s_{j}(p,b),
\end{equation*}
where the last equality follows from $s_j = 1/t_j$. When $p$ and $b$ are real, the last statement follows from $\nabla \log\lvert h\rvert = (\nabla h)/h$.
\end{proof}

Taking the gradient of $(\nabla h(p))/h(p)\cdot b$ with respect to $b$ gives $(\nabla h(p))/h(p)$. Taking the Jacobian with respect to $p$ gives $J_{p}(\nabla h(p)/h(p))$. This proves the following result.

\begin{theorem}\label{thm:gradient_and_hessian}
For $p\in \bC^{k}\smallsetminus \sH$ and a generic direction $b\in \bC^{k}$,
\begin{equation*}
\frac{\nabla h(p)}{h(p)} = -\sum_{j=1}^{\deg h} \nabla_b\, s_{j}(p,b)\qquad\text{and}\qquad J_{p}\Bigl(\frac{\nabla h(p)}{h(p)}\Bigr)= -\sum_{j=1}^{\deg h}J_b\bigl(\nabla_{p} s_{j}(p,b)\bigr).
\end{equation*}
When $p$ is real, the left-hand sides agree with $\nabla \log \lvert h(p)\rvert$ and $\mathrm{Hess}(\log \lvert h(p)\rvert)$, respectively.
\end{theorem}

The first order derivatives of $s_j$ are obtained by solving the linear systems in \eqref{eq:gradient-wrt-p-b}. The second order derivatives are obtained by solving the linear system \eqref{eq:hessian-system}. We summarize the strategy to evaluate $(\nabla h)/h$ and $J_{p}((\nabla h)/h)$ in~\Cref{alg:evaluating_gradient_and_hessian_logh}. 

\medskip
\begin{algorithm}[htbp]
\caption{Evaluating the gradient $(\nabla h)/h$ and its Jacobian $J_{p}((\nabla h)/h)$}
\label{alg:evaluating_gradient_and_hessian_logh}
\SetAlgoLined
\KwIn{A reduced pseudo-witness set $(F,\pi,\mathcal L\times \mathcal L',\mathcal{W})$ for a hypersurface $\sH\subset\bC^k$ as in~\Cref{sec:pseudo-witness}, defined by an unknown nonzero polynomial $h\in\mathbb{R}[p_1,\ldots,p_k]$. A generic point $p\in\bC^k\smallsetminus\sH$.}
\KwOut{The gradient $(\nabla h)/h$ and its Jacobian $J_{p}((\nabla h)/h)$ evaluated at $p$.}
\BlankLine
\textbf{(1) Move $\mathcal{L}\subset \bC^{k}$ through $p$.} Pick a random direction $b\in\bC^k$, move the line $\mathcal{L}$ in the input pseudo-witness set to $\mathcal{L}=\{p+tb\mid t\in\bC\}$, and track the points in $\mathcal{W}$ along this deformation. Let $s_j(p,b)$ and $y_j(p,b)$ for $j=1,\ldots,\deg h$ be as in  \eqref{eq:definition_of_sj_and_yj}.

\BlankLine
\textbf{(2) Solve linear systems.}
For each $j=1,\ldots,\deg h$, compute $\nabla_p s_j$, $J_p y_j$, $\nabla_b s_j$, and $J_b y_j$ by solving the linear system \eqref{eq:gradient-wrt-p-b}, and compute $J_b(\nabla_p s_j)$ by solving \eqref{eq:hessian-system}.

\BlankLine
\textbf{(3) Compute the gradient $(\nabla h)/h$ and its Jacobian $J_{p}((\nabla h)/h)$, and  evaluate them at $p$.} 
Use the formulas from \Cref{thm:gradient_and_hessian}.
\end{algorithm}

\begin{example}
We consider again the quadratic discriminant from \Cref{ex:quadratic_discriminant}. 
As in \Cref{ex:quadratic_discriminant3} we consider a general line $\mathcal L\subset \mathbb C^2$ and $\mathcal L' = \mathbb C$, so that 
$\mathcal{L}\times \mathcal{L}' =  \mathcal{L}\times \bC.$
Let $p = (p_{1},p_{2})\in\bC^2\smallsetminus\sH$ and $b = (b_{1},b_{2}) \in \bC^{2}\smallsetminus \{(0,0)\}$ so that $\mathcal L = \{p+tb\}$. Plugging into \eqref{eq:QuadDisc1} we then have 
    \begin{equation}\label{eq:slice}
    \begin{split}
    G(s,y) & = 
        F\bigl(p + (1/s)b, y\bigr)\\[0.5em] 
        & = 
        \begin{pmatrix}
        y^{2}+(p_{1}+ (1/s) b_{1})y+(p_{2}+(1/s) b_{2})\\
        2y+(p_{1}+(1/s) b_{1})
        \end{pmatrix}
        = \begin{pmatrix}0\\0\end{pmatrix}.
    \end{split}
    \end{equation}
For fixed $(p,b)$, this system has two solutions corresponding to points of $\sH \cap \mathcal{L}$. 

We use the strategy laid out in \Cref{thm:gradient_and_hessian} to compute $J_{p}(\nabla h(p)/h(p))$. We differentiate~\eqref{eq:slice} with respect to $p$ and  $b$ to obtain equations \eqref{eq:gradient-wrt-p-b} with 
Jacobian matrices
    \begin{equation*}
    \begin{split}
    J_{(s,y)}G & =
    \begin{pmatrix}
    (-1/s^{2})(b_1y+b_2) & 2y+p_1+(1/s) b_1\\
    (-1/s^{2})b_1 & 2
    \end{pmatrix},\\[0.5em]
        J_{p}G & = 
        \begin{pmatrix}
            y & 1\\
            1 & 0\\
        \end{pmatrix},\quad  \text{and}
        \quad
        J_{b}G = 
        \begin{pmatrix}
            y/s & 1/s\\
            1/s & 0\\
        \end{pmatrix}.
    \end{split}
    \end{equation*}
    We evaluate these matrices at the zeros of $G$ to determine $\nabla_{p}s_{j}$ and $\nabla_{b}s_{j}$.
    
    Next, computing the systems \eqref{eq:hessian-system} for $i=1,2$ and concatenating them into one system yields
    \begin{equation*}
        J_{(s,y)}G
        \begin{pmatrix}
           \partial_{b_{1}}(\nabla_{p}s) & \partial_{b_{2}}(\nabla_{p}s)\\
           \partial_{b_{1}}(\nabla_{p}y) & \partial_{b_{2}}(\nabla_{p}y)
        \end{pmatrix}
        = -
        \begin{pmatrix}
        R^{(1)}(p,b) & R^{(2)}(p,b)    
        \end{pmatrix}.
    \end{equation*}
    Solving this system yields $J_{b}(\nabla_{p}s) \in \bC^{2\times 2}$.
    Summing over all $s_j$ we obtain, according to the formula in \Cref{thm:gradient_and_hessian}, $J_{p}(\nabla h(p)/h(p))$ $= -J_{b}(\nabla_{p}s_{1}) - J_{b}(\nabla_{p}s_{2})$.\qedhere
\end{example}

Combining~\Cref{thm:gradient_and_hessian} with the definition of $\Phi$ gives the following formulas.

\begin{corollary}
Let $r = \lvert h\rvert/q^{e}$ be a routing function on $\bR^{k}\smallsetminus \sH$ with $q$ as in \eqref{eq:routing-function} and let $\Phi$ be as in~\eqref{eq:Phi}. Then, for $p \in \bC^{k}\smallsetminus (\sH \cup V(q))$, it holds that
\begin{align*}
\Phi(p) &= -\sum_{j=1}^{\deg h} \nabla_b\, s_{j}(p,b)
- \frac{2e(p-c)}{q(p)},\quad  \text{and}\\
J_{p}\Phi(p)
&= -\sum_{j=1}^{\deg h}J_b\bigl(\nabla_{p} s_{j}(p,b)\bigr) + \frac{4e}{q(p)^{2}}(p-c)(p-c)^{T} - \frac{2e}{q(p)}I_{k}.
\end{align*}
When $p$ is real, the left-hand sides agree with $\nabla\log r(p)$ and $\Hess(\log r(p))$, respectively.
\end{corollary}

\subsection{Summary of the algorithm for an unknown defining polynomial}
\label{subsec:summary}

We summarize our strategy for computing complements of hypersurfaces defined by pseudo-witness sets in \Cref{alg:complete}. An implementation of the algorithm in Julia, built on the software package \texttt{HomotopyContinuation.jl} \cite{homotopycontinuation}, will be provided in a forthcoming software paper. The computational bottleneck in \Cref{alg:complete} is solving the system $\Phi(p) = 0$ with monodromy, as each monodromy loop requires many evaluations of $\nabla h/h$ and $J_{p}((\nabla h)/h)$. 

\medskip
\begin{algorithm}[htbp]
\caption{Computing regions of a hypersurface with unknown defining polynomial}
\label{alg:complete}
\SetAlgoLined

\KwIn{A system $F(p,z)\in(\bR[p_1,\ldots,p_k,z_1,\ldots,z_{n-k}])^{m}$ and a union of irreducible components $X\subset V(F)$ such that $\sH = \overline{\pi(X)}$ is a hypersurface, where $\pi(p,z)=p$.}

\KwOut{The regions of the real complement $\bR^k\smallsetminus\sH$.} 
\BlankLine

\BlankLine 

\textbf{(1) Compute a reduced pseudo-witness set}. Find a reduced pseudo-witness set $(F,\pi,\mathcal{L}\times \mathcal L',\mathcal{W})$ of $\mathcal{H}$, and compute the degree $\deg \sH$ as explained in \Cref{sec:pseudo-witness}. 

\BlankLine
\textbf{(2) Construct a routing function.}
Sample $c\in\mathbb{R}^k$ and take $e>\deg(\sH)/2$. This defines a routing function $r$. For the next steps in the algorithm, it is enough to be able to evaluate the rational functions $\Phi(p)$ and $J_{p}\Phi(p)$ defined as in \eqref{eq:Phi}. Use the method in \Cref{subsec:evaluating} for this.

\BlankLine
\textbf{(3) Compute routing points.}
Compute all complex solutions to $\Phi(p) = 0$ through the monodromy approach described in \Cref{sec:HCMonodromy}. The real solutions in $\bR^{k}\smallsetminus \sH$ are the routing points of $r$.

\BlankLine
\textbf{(4) Compute regions.}
Run \Cref{alg:connected-components}.
\end{algorithm}

\section{Case studies}
\label{sec:examples}
We end the paper by using our method to analyze the possible root counts of three parametric polynomial systems arising in applications: the Kuramoto model (\Cref{subsec:kuramoto}), the 3RPR mechanism from kinematics (\Cref{subsec:3RPR}), and
a model for the Allee effect in population ecology (\Cref{subsec:allee}).
Each example consists of a parametric polynomial system 
\begin{equation*}
    G(p;z) \in (\bR[p_{1},\dots,p_{k}, z_{1},\dots,z_{n-k}])^{n-k},
\end{equation*}
with parameters $p=(p_1,\ldots,p_k)$ and variables $z=(z_1,\ldots,z_{n-k})$. We study its \emph{discriminant variety}, which is the hypersurface
\begin{equation}
\label{eq:discriminant}
    \Delta = \overline{\pi\bigl(\bigl\{(p,z) \in \bC^{k}\times \bC^{n-k} \mid G(p;z) = 0,\: \det\bigl(J_zG(p;z)\bigr) = 0\bigr\}\bigr)}
\subset \bC^{k}
\end{equation}
where $\pi\colon \bC^{k}\times \bC^{n-k} \to \bC^{k}$ is the projection to the parameter space.  
Each of the systems considered in this section has the property that no real solution can go to infinity as the parameters approach a finite point in parameter space. Hence, the number of real solutions is finite and constant on each region of the real complement $\bR^{k}\smallsetminus \Delta$. In what follows, we compute the regions using the method outlined in \Cref{alg:complete}. 

\subsection{Kuramoto model}\label{subsec:kuramoto}

The Kuramoto model is an oscillator model that captures synchronization phenomena~\cite{Kuramoto,Kuramoto:1975ebm,Kuramoto-chemical}. It is described by the system of ordinary differential equations
\begin{equation}\label{eq:kuramoto-ode}
    \frac{\mathrm d\theta_{i}}{\mathrm dt} = \omega_{i} - \frac{K}{N}\sum\limits_{j=1}^{N}a_{i,j}\sin(\theta_{i} - \theta_{j}), \qquad \text{for } i=1,\dots,N,
\end{equation}
where $\theta=(\theta_1,\ldots,\theta_N)$ are the phases of the oscillators, $K$ is the coupling strength, $\omega = (\omega_{1},\dots,\omega_{N})$ with $\sum_{i=1}^N\omega_i=0$ are the natural frequencies, and $A=(a_{i,j})$ is the adjacency matrix of an undirected \emph{coupling graph}. 

To study the equilibria modulo rotational symmetry, we fix $\theta_{N} = 0$. Using the identity $\sin(x-y) = \sin(x)\cos(y) - \sin(y)\cos(x)$ we substitute $s_{i} = \sin(\theta_{i})$ and $c_{i} = \cos(\theta_{i})$ to rewrite the remaining $N-1$ equations into polynomial form. Coupling these equations with the Pythagorean identity $s_{i}^{2} + c_{i}^{2} = 1$, we obtain a system of $2(N-1)$ polynomial equations in $2(N-1)$ variables, where $c_N=1$ and $s_N=0$:
\begin{equation}\label{eq:kuramoto-polynomial}
    \begin{split}
        0 & = \omega_{i} - \frac{K}{N}\sum\limits_{j=1}^{N}a_{i,j}(s_{i}c_{j} - s_{j}c_{i}),\\
        0 & = s_{i}^{2} + c_{i}^{2} - 1,
    \end{split}
\end{equation}
for $i=1,\ldots,{N-1}$. The real solutions to~\eqref{eq:kuramoto-polynomial} now correspond to equilibria of the Kuramoto flow. 

We study the Kuramoto model with three oscillators, where the graph is the triangle (i.e., $a_{i,j} = 1$ if $i\neq j$ and $0$ otherwise) and $K=1$. In the formulation of \eqref{eq:kuramoto-polynomial}, this gives rise to the system
\begin{equation}\label{kuramoto_eq}
G(\omega; z) = \begin{pmatrix}
(s_{1}c_{2} - c_{1}s_{2}) + s_{1} - 3 \omega_{1}\\[0.2em]
(s_{2}c_{1} - c_{2}s_{1}) + s_{2} - 3 \omega_{2}\\[0.2em]
s_{1}^{2} + c_{1}^{2} - 1\\[0.2em]
s_{2}^{2} + c_{2}^{2} - 1
\end{pmatrix}
\end{equation}
with parameters $\omega= (\omega_1,\omega_2)$ and variables $z=(c_1,c_2,s_1,s_2)$. 
The discriminant can be computed symbolically, and is given by a degree 12 polynomial in the parameters with 41 terms. Its real complement is known to have 9 connected components with real root counts 0, 2, 4 and 6 \cite{CossHauensteinHongMolzahn2018,RealMonodromy}. 

\begin{figure}[htbp]
\centering
\begin{subfigure}{0.49\textwidth}
\centering
\includegraphics[height=6cm]{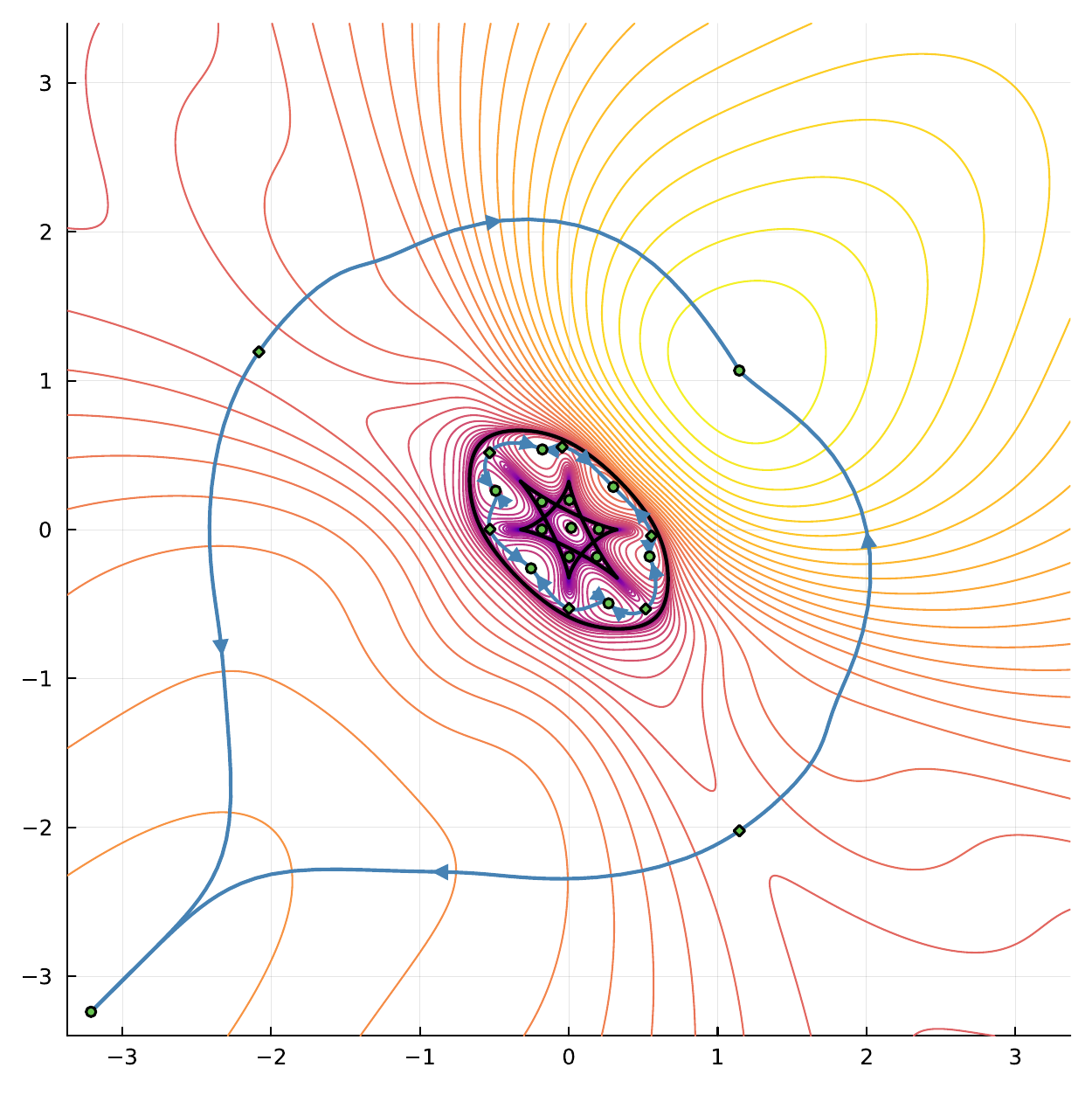}
\end{subfigure}
\begin{subfigure}{0.49\textwidth}
\centering
\includegraphics[height=6cm]
{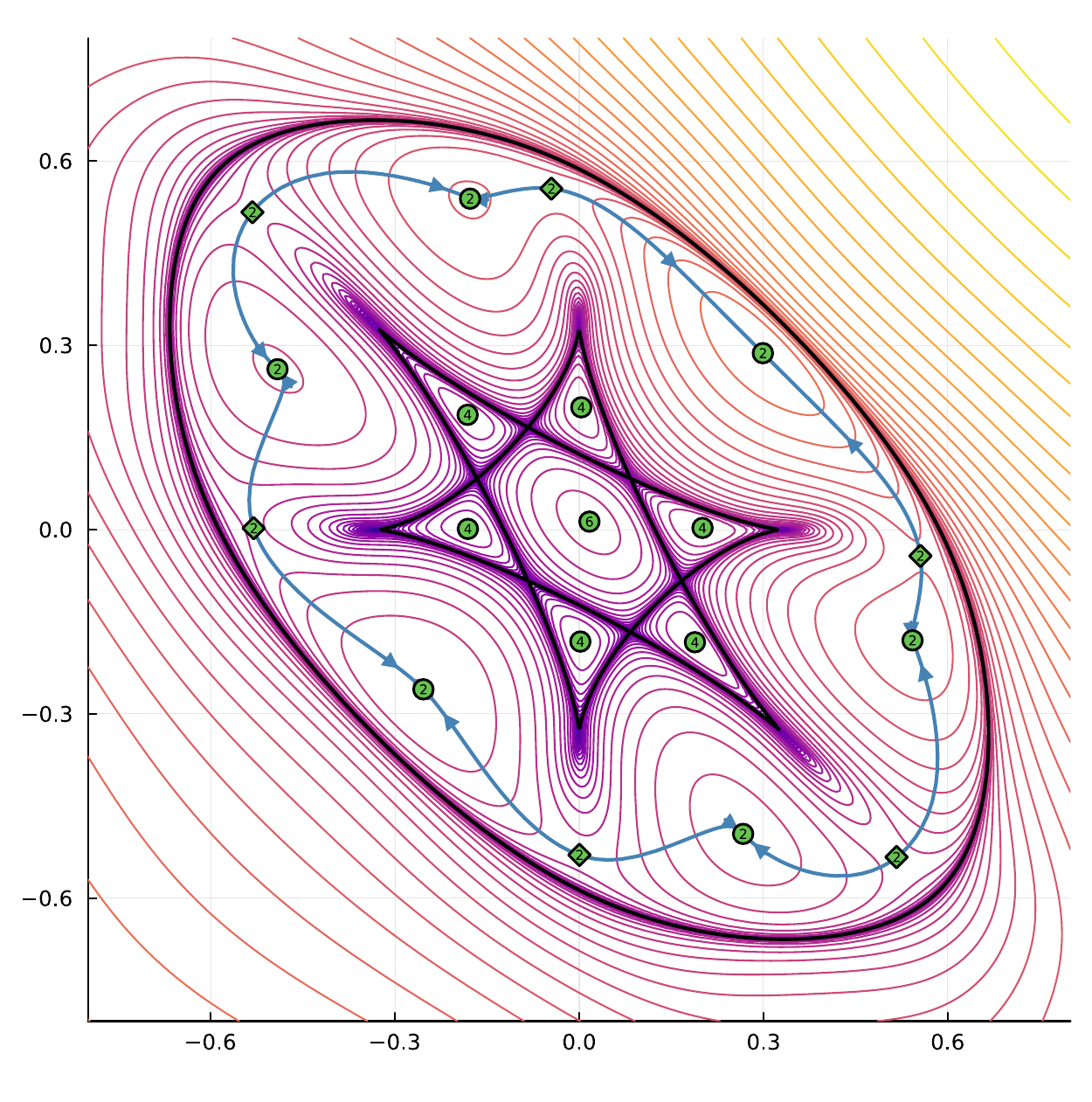}
\end{subfigure}
\caption{Routing points and connecting paths for the discriminant of the Kuramoto model with three oscillators \eqref{kuramoto_eq}. 
The right picture shows a zoomed-in version of the left picture, with the routing points labeled by the corresponding real root counts of \eqref{kuramoto_eq}. The gradient flow is visible in blue, with arrows pointing towards index-zero routing points. We see that $7$ of the $9$ regions are within the star-shaped figure. The remaining two regions are the bounded region enclosed by the oval and the unbounded exterior region. The pictures also display the level sets of the logarithmized routing function $\log r$. The routing points and the flows were computed without access to the discriminant polynomial.}
\label{fig:kuramoto}
\end{figure}

We recover this result via our method, without making use of the discriminant polynomial. Using a routing function with center $(0.47,0.43)$, we find $24$ routing points, and using \Cref{alg:connected-components}, the routing points are grouped into 9 regions. \Cref{fig:kuramoto} shows the routing points and gradient flow connections, along with the known discriminant variety.

\subsection{3RPR}
\label{subsec:3RPR}

The 3RPR mechanism is an important example in kinematics of a system that exhibits a {\em nonsingular assembly-mode change}: the mechanism can pass from one solution branch to another without crossing a singularity~\cite{bonev,GOSSELIN1992107,hayes,husty,Innocenti,macho,RealMonodromy}. It consists of three prismatic legs with revolute joints, each anchored to the ground at one end and attached to a moving triangular platform at the other.

To model its motion, two coordinate frames are used: a fixed frame attached to the ground, where the anchor points lie at fixed positions, and a moving frame attached to the triangular platform, where the platform's corner points remain fixed relative to each other. 
In \Cref{fig:3rpr-mechanism}, the leg lengths are $\ell_{1},\ell_{2}$, and $\ell_{3}$, and we work with their squares $c_{i} = \ell_{i}^2$ for convenience. In the fixed frame, the anchor points are located at $(0, 0)$, $(A_2, 0)$, and $(A_3, B_3)$ while in the moving frame, the platform corners are located at the points $P_{1} = (0,0)$, $P_{2} = (a_{2},0)$, and $P_{3} = (a_{3},b_{3})$.

\begin{figure}[htbp]
\centering
\begin{tikzpicture}[line cap = round, line join = round,scale=0.8]
    \coordinate (A) at (2,2);
    \coordinate (B) at (0,4.5);
    \coordinate (C) at (3,4);

    \coordinate (G1) at (4,1);
    \coordinate (G2) at (0,1);
    \coordinate (G3) at (3,3);

    \draw[line width=1pt] (A) -- (G1);
    \draw[line width=1pt] (B) -- (G2);
    \draw[line width=1pt] (C) -- (G3);

    \fill[black!20, opacity=0.7] (A) -- (B) -- (C) -- cycle;
    \draw[line width=1pt, black, opacity=0.7] (A) -- (B) -- (C) -- cycle;

    \fill[black] (A) circle(2pt);
    \fill[black] (B) circle(2pt);
    \fill[black] (C) circle(2pt);

    \def\t{0.12}
    \fill (G1) ++(-\t,-\t) -- ++(2*\t,0) -- ++(-\t,2*\t) -- cycle;
    \fill (G2) ++(-\t,-\t) -- ++(2*\t,0) -- ++(-\t,2*\t) -- cycle;
    \fill (G3) ++(-\t,-\t) -- ++(2*\t,0) -- ++(-\t,2*\t) -- cycle;

    \put(-15,60){$\ell_{1}$}
    \put(70,40){$\ell_{2}$}
    \put(73,75){$\ell_{3}$}
    \put(-32,22){$(0,0)$}
    \put(100,22){$(A_{2},0)$}
    \put(75,60){$(A_{3},B_{3})$}    
\end{tikzpicture}
\caption{An example of a 3RPR mechanism.}
\label{fig:3rpr-mechanism}
\end{figure}

\begin{figure}[htbp]
\centering
\begin{subfigure}{0.48\textwidth}
\centering
\includegraphics[height=5.5cm]{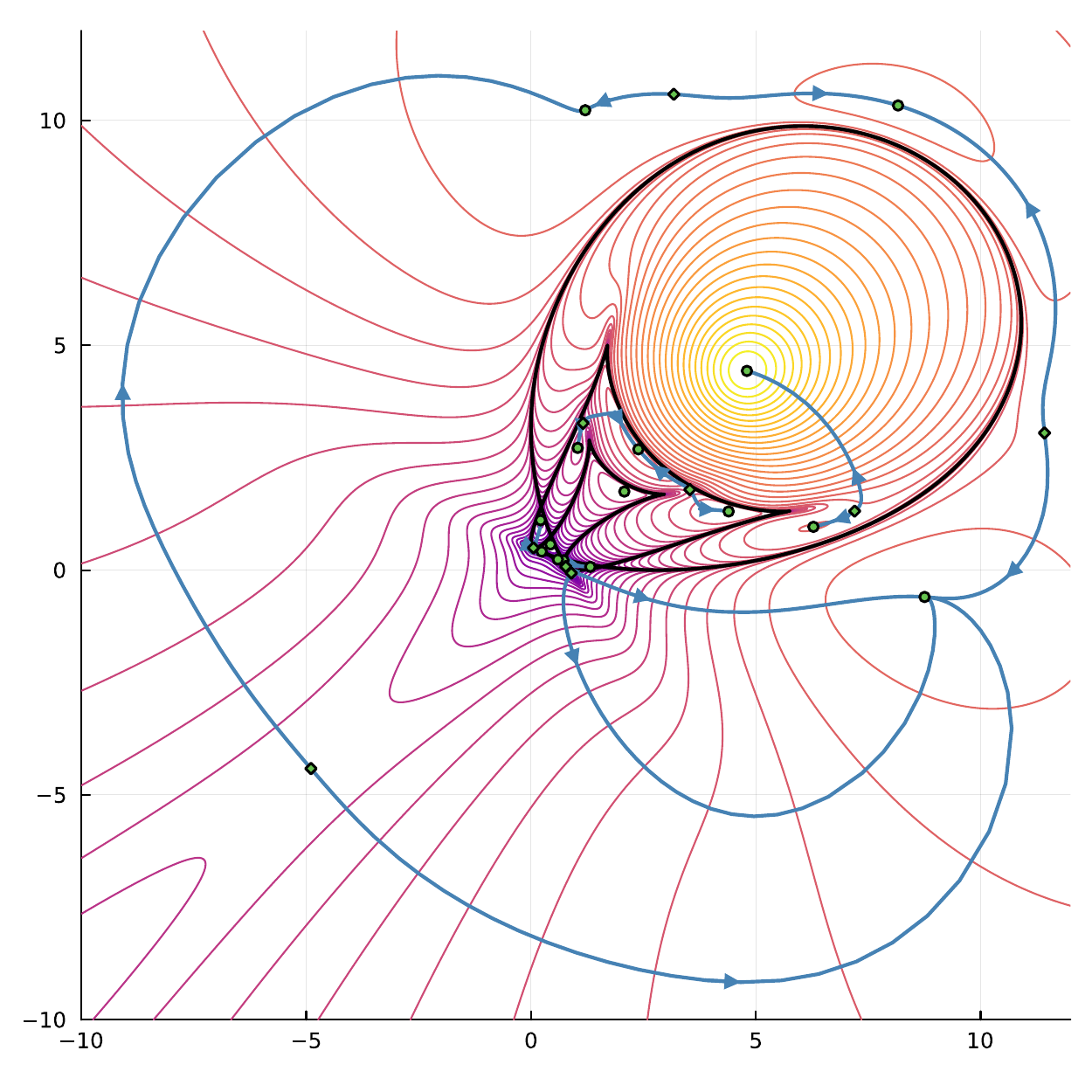}
\end{subfigure}
\begin{subfigure}{0.48\textwidth}
\centering
\includegraphics[height=5.5cm]
{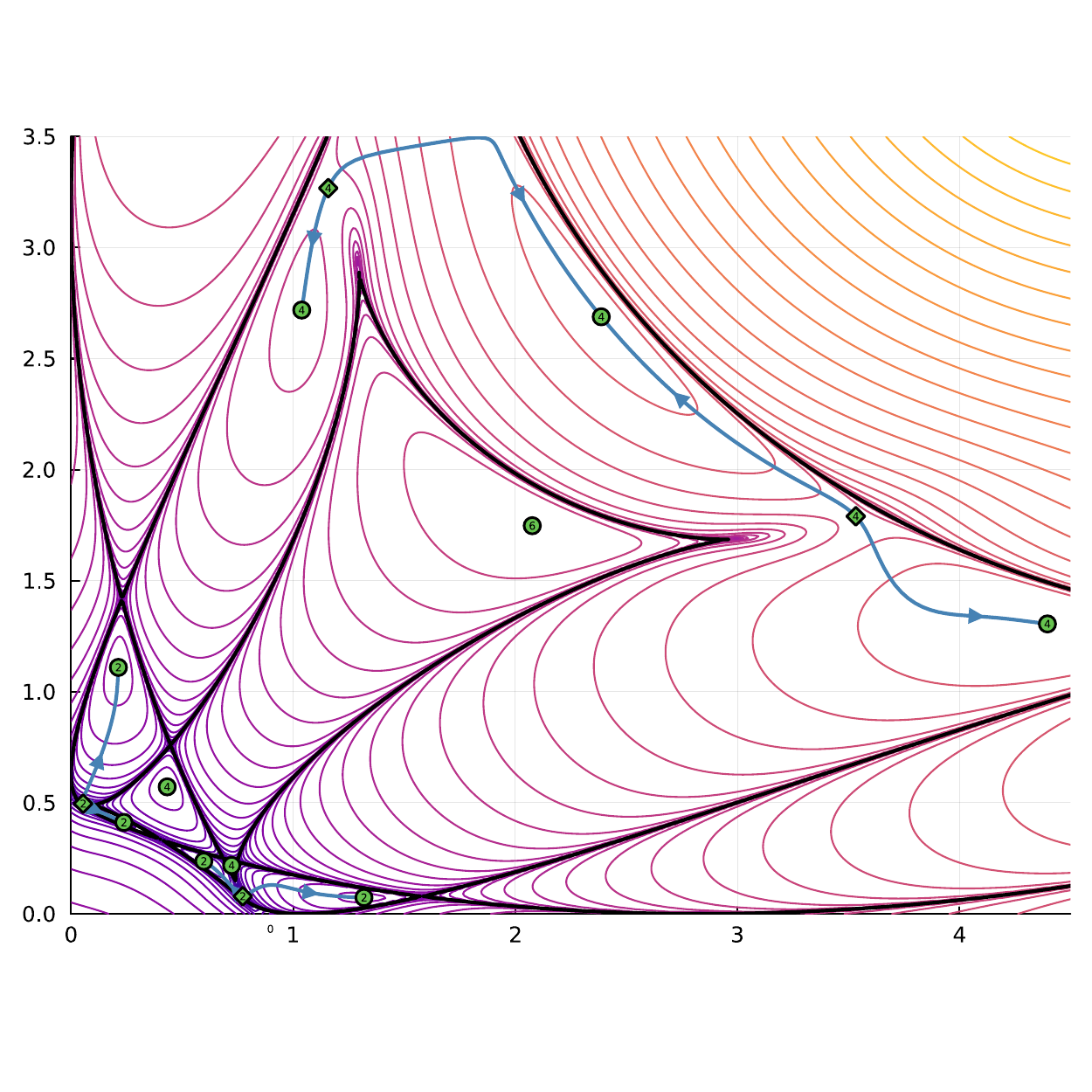}
\end{subfigure}
\caption{Routing points and connecting paths for the two-parameter version of the 3RPR system.  The right picture shows a zoomed-in version of the left picture, with the routing points labeled by the corresponding real root counts. The gradient flow is visible in blue. The arrows point towards index-zero routing points. The pictures also display the level sets of the logarithmized routing function $\log r$. The routing points and the flows were computed without access to the discriminant polynomial.}
\label{fig:3RPRv0}
\end{figure}

A fundamental problem is to determine all the ways a real motion of the mechanism can start and end in the same ``home'' configuration for a fixed set of leg lengths. The configuration variables are $(p,\phi) = (p_{1},p_{2},\phi_{1},\phi_{2})$, where $p = (p_{1},p_{2})$ gives the coordinates of the platform point $P_{3}$ in the fixed frame and $\phi = (\phi_{1},\phi_{2})$ is a unit circle representation of the platform's rotation. This yields a polynomial system $G(a,A,b_{3},B_3,c;p,\phi)$ with variables $(p,\phi)$, parameters 
\begin{equation*}
a = (a_{2},a_{3}), \quad A=(A_{2},A_{3}), \quad b_{3}, \quad B_{3}, \quad \text{and} \quad c=(c_{1},c_{2},c_{3}),
\end{equation*}
and polynomials
\begin{align*}
g_1 & = \phi_1^2 + \phi_2^2 - 1\\
g_2 & = p_1^2 + p_2^2 - 2(a_3p_1 + b_3p_2)\phi_1 + 2(b_3p_1 - a_3p_2)\phi_2 + a_3^2 + b_3^2 - c_1\\
g_3 & = p_1^2 + p_2^2 - 2A_2p_1 +
2((a_2-a_3)p_1 - b_3p_2 + A_2a_3 - A_2a_2)\phi_1 \\
&\qquad +~2(b_3p_1 + (a_2-a_3)p_2 - A_2b_3)\phi_2
+ (a_2-a_3)^2 + b_3^2 + A_2^2 - c_2\\
g_4 &= p_1^2 + p_2^2 - 2(A_3 p_1 + B_3 p_2) + A_3^2 + B_3^2 - c_3.
\end{align*}
Following a scaled version of the formulation in~\cite{husty} and \cite[Section~4]{RealMonodromy}, we consider a two-parameter version of the problem obtained by viewing $(c_1,c_2)$ as parameters and fixing
\begin{equation}
\label{eq:3RPR_fixed_parmaters}
a_2 = 1.4,\quad a_3 = 0.7,\quad  b_3 = 1.0,\quad  A_2 = 1.6,\quad A_3 = 0.9,\quad B_3 = 0.6,\quad \text{and } c_3 = 1.0\,.
\end{equation}
In this case, it is known that the discriminant is a degree 12 curve, and that its real complement has 8 connected components with real root counts 0, 2, 4 and 6 \cite[Section~4]{RealMonodromy}. We find all of them with our method, as illustrated in \Cref{fig:3RPRv0} (which can be compared with \cite[Fig.~8]{RealMonodromy}). More precisely, using a routing function with center $(4.72, 4.33)$, we find 24 routing points, of which 15 have index zero.

Next, we consider a three-parameter version of the problem, where we view $(c_1,c_2,c_3)$ as free parameters, and fix the remaining ones according to \eqref{eq:3RPR_fixed_parmaters}. In this case,  symbolic elimination shows that the discriminant surface is given by a polynomial of degree 12 with 455 terms. Using a routing function centered at $(4.72, 4.33, 1.70)$, our method finds 60 routing points and 9 regions of the complement of the real part of the discriminant. In these regions, the real root counts 0, 2, 4 and 6 are recorded. This is in agreement with the result obtained by analyzing the symbolic polynomial with \texttt{HypersurfaceRegions.jl} \cite{HypersurfaceRegions}.

Finally, we consider another three-parameter version of the problem, where we instead view $(c_1,c_2,A_2)$ 
as free parameters, and fix other parameters as in \eqref{eq:3RPR_fixed_parmaters}. In this case, we have not been able to find a symbolic expression for the discriminant surface, but its degree can readily be determined to be 24 through a pseudo-witness set. 

This example is at the limit of what our current implementation can handle, as the monodromy step did not terminate after 14 days\footnote{Computations were run on a 4-socket Intel Xeon Gold 6128 system (24 cores total, 3.4 GHz) with 3 TB RAM.}. At this point, 457 complex critical points for the routing function centered at $( 4.72, 4.33, 1.70)$ had been found, of which 90 points were real. While likely not being a complete set of routing points, these points still provide a valuable sample of the parameter space. In particular, we can use them to verify that there are open regions in the real parameter space with real root counts of 0, 2, 4 and 6, respectively.

\subsection{Allee effect}\label{subsec:allee}
In population ecology, the \emph{Allee effect} refers to the phenomenon in which a population can experience a negative growth rate if its size falls below a critical threshold, unlike classical logistic growth. There are many possible factors that can cause this effect, and it has been widely studied (see, e.g., \cite{berec2007multiple}). We consider the polynomial model from \cite{alleethreepopulations}, where the Allee effect for three \emph{patches} $z=(z_1,z_2,z_3)$ is modeled by the dynamical system 
\begin{align*}
\frac{\mathrm d z}{\mathrm d t} = G(a,b,z)\quad \text{ where }\quad  G(a,b,z) = \begin{pmatrix}
    z_1(1-z_1)(z_1-b)+a(z_2-z_1) + a(z_3 - z_1)\\[0.25em]
    z_2(1-z_2)(z_2-b)+a(z_1-z_2) + a(z_3 - z_2)\\[0.25em]
    z_3(1-z_3)(z_3-b)+a(z_1-z_3) + a(z_2 - z_3)\end{pmatrix}.
\end{align*}
As discussed in \cite{alleethreepopulations}, a \emph{patch} can be interpreted as a habitat, a cellular compartment,
or a microplate. The parameter $a > 0$ models the \emph{dispersal rate} between patches, and $0<b < \frac{1}{2}$ is the \emph{Allee
threshold}.

The steady states of this dynamical system are the solutions of the system of equations 
$$G(a,b,z)  = 0$$
with parameters $p=(a,b)$ and variables $z=(z_1,z_2,z_3)$. We are interested in exploring the possible number of \emph{nonnegative} solutions in the \emph{biologically relevant} parameter space $(0,\infty)\times (0,\frac{1}{2})$. 

To this end, we extend the discriminant $\Delta$ from \eqref{eq:discriminant} to $\Delta\cup V(ab(b-\frac{1}{2}))$. As explained in \Cref{extend_h},
on the level of routing functions, this means adding $\log(\vert a\vert)$,  $\log(\vert b\vert)$ and $\log(\vert b-\frac{1}{2}\vert)$ to $\log(r)$. It is clear from inspection of $G$ that $z=0$ is always a solution for biologically relevant parameters, and that all other nonnegative solutions have strictly positive coordinates. Hence, the number of nonnegative solutions is constant in each region of $((0,\infty)\times (0,\frac{1}{2}))\smallsetminus (\Delta\cup V(ab(b-\frac{1}{2})))$.

Computing a defining equation for $\Delta$ is highly challenging~\cite{sadeghimanesh2022resultanttoolsparametricpolynomial,song2025steadystateclassificationallee,tsai2026}, but the relevant regions have nevertheless been determined through cylindrical algebraic decomposition in \cite{alleethreepopulations}. We demonstrate that the regions can also be computed using our method.

\begin{figure}[htbp]
\centering
\begin{subfigure}{0.48\textwidth}
\centering
\includegraphics[width = 0.9\textwidth]{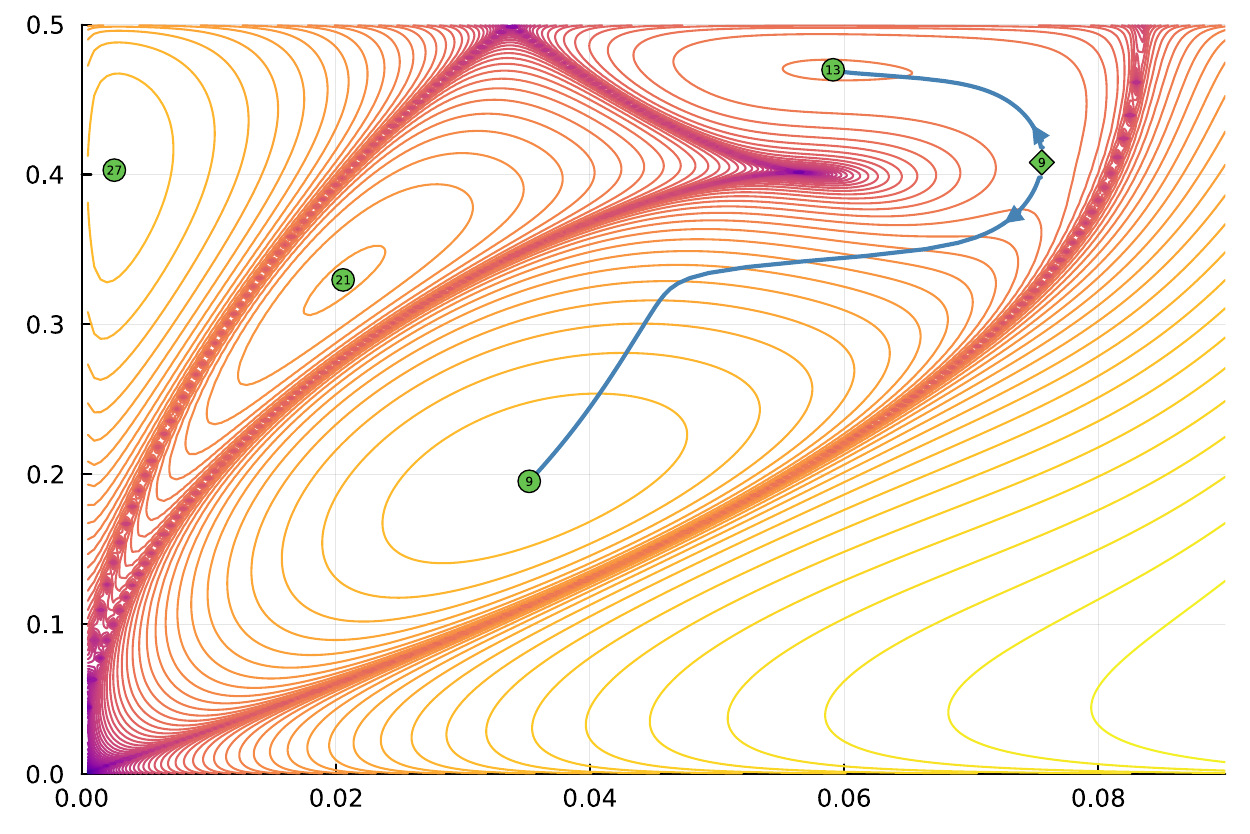}
\end{subfigure}
\begin{subfigure}{0.48\textwidth}
\centering
\includegraphics[width = 0.9\textwidth]
{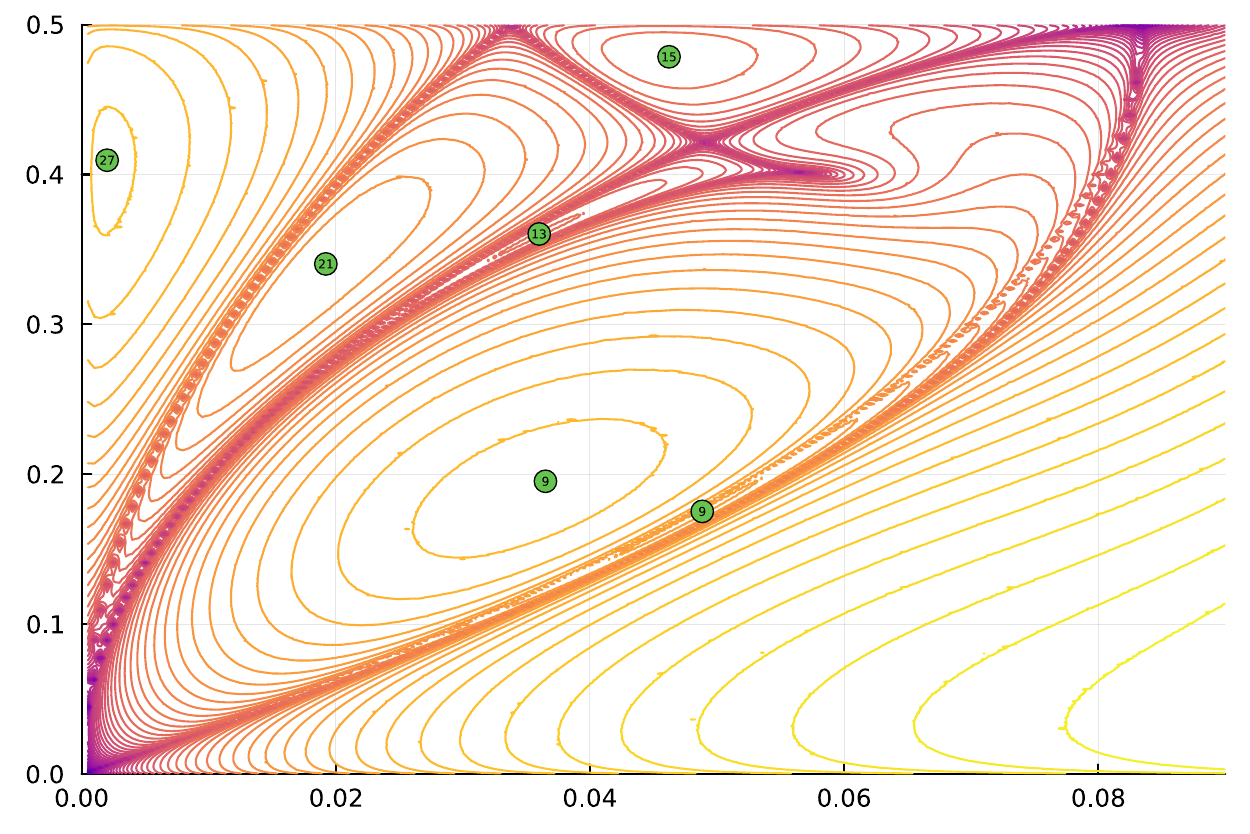}
\end{subfigure}
\caption{Contours and critical points of a routing function for the modified discriminant $\Delta\cup V(ab(b-\frac{1}{2}))$ of the Allee system from \Cref{subsec:allee}. The left figure corresponds to using the formulation \eqref{eq:discriminant}, which gives the subvariety~$\hat{\Delta}$. The right figure uses the radical and gives the full discriminant $\Delta$. The routing points are in green and are labeled by the number of nonnegative solutions of $G(a,b,z)=0$. In both figures, there is a routing point of the unbounded region which is outside the picture, with nonnegative root count 3.  Contours (up to a constant), critical points and gradient flow were computed without access to the discriminant polynomial.}
\label{fig:allee}
\end{figure}

An interesting property of this example is that $\Delta$ is not irreducible, and that some of its irreducible components have a non-reduced pseudo-witness set with respect to the equations in ~\eqref{eq:discriminant}.
Hence, if these equations are used, our implementation effectively works with the subvariety $\hat{\Delta}\subset \Delta$ corresponding to the reduced part. We can see the effect of this in the left part of \Cref{fig:allee}, which shows the resulting routing function and the biologically relevant routing points, labeled by the number of nonnegative solutions of~$G(a,b,z)=0$. As we can see, gradient flow connects a parameter pair $(a,b)$ for which $G(a,b,z)=0$ has 15 nonnegative solutions with another pair, for which there are 9 nonnegative solutions. 

To recover the missing part of the discriminant, we could, in principle, use deflation techniques. In this case, however, it turns out that the radical $\sqrt{I}$ of the ideal 
$$I=\left\langle G(a,b,z),\ \det (J_zG(a,b,z))\right\rangle$$
can readily be computed symbolically in the computer algebra system \texttt{OSCAR} \cite{OSCAR,OSCAR-book}. By using a generating set of $\sqrt{I}$, together with the standard ``squaring up'' technique \cite[Section~2.5]{bates2024numericalnonlinearalgebra} for dealing with overdetermined systems, we obtain a reduced pseudo-witness set for $\Delta$. The resulting routing function and the biologically relevant routing points are shown in the right part of \Cref{fig:allee}, which agrees with the regions displayed in \cite[Figure~3]{alleethreepopulations}.

\section{Conclusion}\label{sec:Conclusion}

Since computing a defining
polynomial for a hypersurface
arising as the projection
of another variety is a 
computationally challenging
elimination problem,
our algorithm (\Cref{alg:complete})
provides a novel method to recover the regions of the real complement of a hypersurface without computing its defining polynomial. This allows one to compute the regions of the complement
even when the hypersurface arises as the projection of an algebraic set. By intersecting this hypersurface with a line, we obtain the needed numerical information to effectively perform ``elimination without eliminating''
as illustrated with several examples. 

\section*{Acknowledgments}\label{sec:Ack}
We thank the Fields Institute for Research in Mathematical Sciences and the organizers of the Workshop on the Applications of Commutative Algebra where the initial steps of this research were conducted. 
We also thank Hannah Friedman for helpful discussions in the early phases of~the~project. 

\section*{Funding}
PB was supported by DFG, German Research Foundation -- Projektnummer 445466444.
JC acknowledges the support of the National Science Foundation Grant DMS-2402199.
AKE was partially supported by National Science Foundation Grant DMS-2023239.
JDH was partially supported
by the 
National Science Foundation Grant CCF-2331400, Simons
Foundation SFM-00005696,
and the Robert and Sara Lumpkins Collegiate Professorship.
OH was partially funded by
the European Union (Grant Agreement no. 101044561, POSALG). Views and opinions expressed
are those of the authors only and do not necessarily reflect those of the European Union or European
Research Council (ERC). Neither the European Union nor ERC can be held responsible for them. DKJ and DM were partially supported by NSERC Discovery Grant RGPIN-2023-03551. 
JLG was partially supported by the Robert and Sara Lumpkins Collegiate Professorship.

\newlength{\bibitemsep}\setlength{\bibitemsep}{.05\baselineskip plus .05\baselineskip minus .05\baselineskip}
\newlength{\bibparskip}\setlength{\bibparskip}{0pt}
\let\oldthebibliography\thebibliography
\renewcommand\thebibliography[1]{
  \oldthebibliography{#1}
  \setlength{\parskip}{6\bibitemsep}
  \setlength{\itemsep}{\bibparskip}
}

\bibliographystyle{abbrv}
\bibliography{ref}

\bigskip

{\samepage

\noindent {\bf Authors' addresses:}\\
Paul Breiding, University of Osnabr\"uck \hfill{\tt pbreiding@uni-osnabrueck.de}\\
John Cobb, Auburn University \hfill{\tt john.cobb@auburn.edu}\\
Aviva K. Englander, University of Wisconsin--Madison \hfill{\tt akenglander@wisc.edu}\\
Nayda Farnsworth, Colgate University \hfill{\tt nfarnsworth@colgate.edu}\\
Jonathan D. Hauenstein, University of Notre Dame \hfill{\tt hauenstein@nd.edu}\\
Oskar Henriksson, MPI-CBG Dresden \hfill{\tt oskar.henriksson@mpi-cbg.de}\\
David K. Johnson, University of Western Ontario \hfill{\tt djohn225@uwo.ca}\\
Jordy Lopez Garcia, University of Notre Dame \hfill{\tt jlopezga@nd.edu}\\
Deepak Mundayur, University of Western Ontario \hfill{\tt dmundayu@uwo.ca}

}

\end{document}